\documentclass{amsart}

\usepackage{amsmath} 
\usepackage{amssymb}
\usepackage{mathrsfs}

\newbox\stickbox
\setbox1=\hbox{\raise1ex\hbox{$\bullet$}}
\setbox0\hbox to \wd1{\hss\vrule width 0.4pt\strut\hss}
\setbox\stickbox\hbox{\box0\box1}
\def\stick{\leavevmode\copy\stickbox}

\newtheorem{theorem}{Theorem}[section] 
\newtheorem{claim}[theorem]{Claim}

\theoremstyle{definition}
\newtheorem{definition}[theorem]{Definition}

\newtheorem{observation}[theorem]{Observation} 

\newtheorem{convention}[theorem]{Convention}

\theoremstyle{remark}
\newtheorem{remark}[theorem]{Remark}
\newtheorem{notation}[theorem]{Notation}

\newcommand{\rest}{{\restriction}}

\newcommand{\wilog}{{\rm without loss of generality}}

\newcommand{\then}{{\underline{then}}}
\newcommand{\when}{{\underline{when}}}
\newcommand{\Then}{{\underline{Then}}}
\newcommand{\If}{{\underline{if}}}
\newcommand{\Iff}{{\underline{iff}}}
\newcommand{\mn}{{\medskip\noindent}}
\newcommand{\sn}{{\smallskip\noindent}}
\newcommand{\bn}{{\bigskip\noindent}}

\newcommand{\cA}{{\mathcal A}}

\newcommand{\gH}{{\mathfrak H}}

\newcommand{\cH}{{\mathcal H}}

\newcommand{\cF}{{\mathcal F}}

\newcommand{\bbP}{{\mathbb P}}
\newcommand{\cP}{{\mathcal P}}

\newcommand{\bbQ}{{\mathbb Q}}

\newcommand{\cU}{{\mathcal U}}

\newcount\skewfactor
\def\mathunderaccent#1#2 {\let\theaccent#1\skewfactor#2
\mathpalette\putaccentunder}
\def\putaccentunder#1#2{\oalign{$#1#2$\crcr\hidewidth
\vbox to.2ex{\hbox{$#1\skew\skewfactor\theaccent{}$}\vss}\hidewidth}}
\def\name{\mathunderaccent\tilde-3 }

\newenvironment{PROOF}[2][\proofname.]
   {\begin{proof}[#1]\renewcommand{\qedsymbol}{\textsquare$_{\rm #2}$}}
   {\end{proof}}

\begin{document}


\title {Universality among graphs omitting a complete bipartite graph} 
\author {Saharon Shelah}
\address{Einstein Institute of Mathematics\\
Edmond J. Safra Campus, Givat Ram\\
The Hebrew University of Jerusalem\\
Jerusalem, 91904, Israel\\
 and \\
 Department of Mathematics\\
 Hill Center - Busch Campus \\ 
 Rutgers, The State University of New Jersey \\
 110 Frelinghuysen Road \\
 Piscataway, NJ 08854-8019 USA}
\email{shelah@math.huji.ac.il}
\urladdr{http://shelah.logic.at}
\thanks{The author thanks Alice Leonhardt for the beautiful typing.
This research was supported by the United States-Israel Binational 
Science Foundation. Publication 706.} 

\date{May 5, 2010}

\begin{abstract}
For cardinals $\lambda,\kappa,\theta$ we consider the class of
graphs of cardinality $\lambda$ which has no subgraph which is
$(\kappa,\theta)$-complete bipartite graph.  The question is whether in such
a class there is a universal one under (weak) embedding.  We solve this
problem completely under GCH.  Under various assumptions mostly related to
cardinal arithmetic we prove non-existence of universals for this
problem.  We also look at combinatorial properties useful for those
problems concerning $\kappa$-dense families.
\end{abstract}

\maketitle
\numberwithin{equation}{section}
\setcounter{section}{-1}

%\newpage

\section* {Anotated Content}
\bigskip

\noindent
\S0 \quad Introduction
\bigskip

\noindent
\S1 \quad Some no we-universal
\begin{enumerate}
\item[${{}}$]  [We define Pr$(\lambda,\kappa)$ and using it gives
sufficient conditions for the non-existence of we-universal in 
${\gH}_{\lambda,\theta,\kappa}$ mainly when $\theta = \kappa^+$.  Also
we give sufficient conditions for no ste-universal when $\lambda =
\lambda^\kappa,2^\kappa \ge \theta \ge \kappa$.]
\end{enumerate}
\bn
\S2 \quad No we-universal by Pr$(\lambda,\kappa)$ and its relatives
\begin{enumerate}
\item[${{}}$]  [We give finer sufficient conditions and deal/analyze
the combinatorial properties we use; Pr$(\lambda,\kappa)$ says that
there are partial functions $f_\alpha$ from $\lambda$ to $\lambda$ for
$\alpha < \lambda$, which are dense, $\kappa = \text{\rm otp}(\text{\rm
Dom}(f_\alpha))$ and $\kappa > 
\text{\rm otp}(\text{\rm Dom}(f_\alpha) \cap 
\text{ Dom}(f_\beta))$ for $\alpha \ne \beta$.]
\end{enumerate}
\bn
\S3 \quad Complete characterization under G.C.H.
\begin{enumerate}
\item[${{}}$]  [Two theorems cover G.C.H. - one assume $\lambda$ is
strong limit and the other assume $\lambda = \mu^+ = 2^\mu + 2^{< \mu}
= \mu$.  Toward this we prove mainly some results on existence of
universe.]
\end{enumerate}
\bn
\S4 \quad More accurate properties

\section {Introduction} 

On the problem of ``among graphs with $\lambda$ nodes and no complete
subgraph with $\kappa$ nodes, is there a universal one'' (i.e. under weak
embedding) is to a large extent solved in Komjath-Shelah
\cite{KoSh:492}, see more there.  E.g. give a complete solution
under the assumption of GCH.

Now there are some variants, mainly for
graph theorists embedding, i.e. a one to one function mapping an edge to
an edge, called here weak or we-embedding; for model theorist an
embedding also maps a non-edge to a non-edge, call strong or ste-embedding.
We have the corresponding we-universal and ste-universal.
\bn 
We deal here with the problem ``among the graphs with $\lambda$ nodes
and no complete $(\theta,\kappa)$-bipartite sub-graph, is there a
universal one?", see below on earlier results. 
We call the family of such graphs ${\gH}_{\lambda,\theta,\kappa}$,
and consider both the weak embedding (as most graph theorists use) and
the strong embedding.  Our neatest result appears in section 3 (see \ref{2.9},
\ref{2.14}).

\begin{theorem}
\label{0.A}  Assume $\lambda \ge \theta \ge \kappa
\ge \aleph_0$. 

\noindent
1) If $\lambda$ is strong limit then: there is a member of 
${\gH}_{\lambda,\theta,\kappa}$ which is we-universal (= universal under weak
embedding) \Iff \, there is a member of 
${\gH}_{\lambda,\theta,\kappa}$ which is
ste-universal (= universal under strong embeddings) \Iff \, {\rm cf}$(\lambda) 
\le { \text{\rm cf\/}}(\kappa)$ and $(\kappa < \theta \vee 
{ \text{\rm cf\/}}(\lambda) < { \text{\rm cf\/}}(\theta))$. 

\noindent
2) If $\lambda = 2^\mu = \mu^+$ and $\mu = 2^{< \mu}$, \then \,
\mn
\begin{enumerate}
\item[$(a)$]  there is no ste-universal in ${\gH}_{\lambda,\theta,\kappa}$
\sn
\item[$(b)$]  there is we-universal in 
${\gH}_{\lambda,\theta,\kappa}$ iff $\mu = \mu^\kappa$ and $\theta =
\lambda$.
\end{enumerate}
\end{theorem}
\mn
We give many sufficient conditions for the non-existence of universals
(mainly we-universal) and some for the existence, for this dealing
with some set-theoretic properties.  Mostly when we get ``no $G \in
{\gH}_{\lambda,\theta,\kappa}$ is we/ste-universal" we, moreover,
get ``no $G \in {\gH}_{\lambda,\theta,\kappa}$ is we/ste-universal
among the bipartite ones".  Hence we get also results on families of
bi-partite graphs.  We do not look at the case $\kappa < \aleph_0$
here.

Rado has proved that: if $\lambda$ is regular $> \aleph_0$ and $2^{<
\lambda} = \lambda$, \then \, ${\gH}_{\lambda,\lambda,1}$ has a
ste-universal member (a sufficient condition for $G^*$ being
ste-universal for ${\gH}_{\lambda,\lambda,1}$ is:  
for any connected graph $G$ with $<\lambda$ nodes, 
$\lambda$ of the components of $G$ are isomorphic to
$G$). Note that $G \in {\gH}_{\lambda,\lambda,1}$ iff $G$ has
$\lambda$ nodes and the valency of every node is $< \lambda$.  Erd\"os
and Rado (see \cite{EH1}, in Problem 74) ask what occurs, under GCH
to say $\aleph_\omega$.  By \cite[3.1]{Sh:26} if $\lambda$ is strong
limit singular then there is a ste-universal graph in 
${\gH}_{\lambda,\lambda,1}$.

Komjath and Pach \cite{KoPa84} prove that
$\diamondsuit_{\omega_1} \Rightarrow$ no universal in 
${\gH}_{\aleph_1,\aleph_1,\aleph_0}$, this holds also for
${\gH}_{\kappa^+,\kappa^+,\kappa}$ when 
$\diamondsuit_{S^{\kappa^+}_\kappa}$ holds; subsequently the 
author showed that $2^\kappa = \kappa^+$ suffice (Theorem 1 there).  
Then Shafir (see \cite[Th.1]{Sha01}) presents this and proves the following:
\bigskip

\noindent
(\cite[Th.2]{Sha01}):  if $\kappa = \text{
cf}(\kappa),\clubsuit_{S^{\kappa^+}_\kappa}$ and there is a MAD family
on $[\kappa]^\kappa$ of cardinality $\kappa$, \then \,
${\gH}_{\kappa^+,\kappa,\kappa}$ has no we-universal.
\smallskip

\noindent
(\cite[Th.3]{Sha01}): if $\kappa \le \theta \le 2^\kappa$ and there is
${\cA} \subseteq [\kappa]^\kappa$ of cardinality $\kappa$ such that
no $B \in [\kappa]^\kappa$ is included in $\theta$ of them, then
${\frak H}_{\theta^\kappa,\kappa,\theta}$ has no ste-universal member.
\smallskip

\noindent
(\cite[Th.4]{Sha01}):  if $\kappa \le \theta \le 2^\kappa$ and
$\clubsuit_{S^\lambda_\kappa}$ then 
${\gH}_{\lambda,\theta,\kappa}$ has no ste-universal members. 
\medskip

\noindent
Here we characterize ``${\frak H}_{\lambda,\theta,\kappa}$ has universal"
under GCH (for weak and for strong embeddings).  We also in \ref{1.1}
prove $\clubsuit_{S^{\kappa^+}_\kappa} \Rightarrow$ no we-universal in
${\gH}_{\kappa^+,\kappa^+,\kappa}$ (compared to \cite[Th.2]{Sha01}, we
omit his additional assumption ``no MAD ${\cA} \subseteq
[\kappa]^\kappa,|{\cA}| = \kappa$"); in \ref{1.2} 
we prove more.  Also (\ref{1.4}) 
$\lambda = \lambda^\kappa \ge 2^\kappa \ge \theta \ge \kappa
\Rightarrow$ no universal under strong embedding in ${\gH}_{\lambda,
\kappa,\theta}$ (compared to \cite[Th.3]{Sha01} we omit an assumption).

Lately some results for which we originally used 
\cite{Sh:460} now instead use \cite{Sh:829}, \cite{Sh:922} which
gives stronger results.
\bigskip

\noindent
\centerline {$* \qquad * \qquad *$}
\bigskip

\noindent
\begin{notation}
\mn
\begin{enumerate}
\item[$\bullet$]   We use $\lambda,\mu,\kappa,\chi,\theta$ for 
cardinals (infinite if not said otherwise)
\sn
\item[$\bullet$]   We use $\alpha,\beta,\gamma,\varepsilon,
\zeta,\xi,i,j$ for ordinals, $\delta$ for limit ordinals
\sn
\item[$\bullet$]   For $\kappa = \text{ cf}(\kappa) <
\lambda,S^\lambda_\kappa = \{\delta < \lambda:\text{ cf}(\delta) = \kappa\}$
\sn
\item[$\bullet$]   $[A]^\kappa = \{B \subseteq A:|B| = \kappa\}$.
\sn
$\cdot$ We use $G$ for graphs and for bipartite graphs; see below Definition
\ref{0.1}(1), \ref{0.2}(1), it will always be clear from the
context which case we intend.
\end{enumerate}
\end{notation}

\begin{definition}
\label{0.1}  
1) A graph $G$ is a pair $(V,R) = (V^G,R^G),V$ 
a non-empty set, $R$ a symmetric irreflexive 2-place relation on it.
We call $V$ the set of nodes of $G$ and $|V|$ is the cardinality of
$G$, denoted by $\|G\|$, and may write 
$\alpha \in G$ instead of $\alpha \in V^G$.

Let $E^G = \{\{\alpha,\beta\}:\alpha R^G \beta\}$, so we may consider $G$ as
$(V^G,E^G)$. 

\noindent
2) We say $f$ is a strong embedding of $G_1$ into $G_2$ (graphs) \If \,:
\mn
\begin{enumerate}
\item[$(a)$]   $f$ is a one-to-one function from $G_1$ into $G_2$;
pedantically from $V^{G_1}$ into $V^{G_2}$
\sn
\item[$(b)_{\text{st}}$]   for $\alpha,\beta \in G_1$ we have 

$\alpha R^{G_1} \beta \Leftrightarrow f(\alpha) R^{G_2} f(\beta)$.
\end{enumerate}
\mn
3) we say $f$ is a weak embedding of $G_1$ into $G_2$ if
\mn
\begin{enumerate}
\item[$(a)$]   above and
\sn
\item[$(b)_{\text{we}}$]   for $\alpha,\beta \in G_1$ we have 

$\alpha R^{G_1} \beta \Rightarrow f(\alpha) R^{G_2} f(\beta)$.
\end{enumerate}
\mn
4)  The $\lambda$-complete graph $K_\lambda$ is the graph
$(\lambda,R)$ were $\alpha R \beta \Leftrightarrow \alpha 
\ne \beta$ or any graph isomorphic to it.
\end{definition}

\begin{definition}
\label{0.2}  
1) $G$ is a bipartite graph means $G =
(U,V,R) = (U^G,V^G,R^G)$ where $U,V$ are disjoint 
non-empty sets, $R \subseteq
U \times V$.  For a bipartite graph $G$, we would like sometimes to
treat as a usual graph (not bipartite), so let $G$ as a graph, 
$G^{[\text{gr}]}$, be $(U^G \cup V^G,\{(\alpha,\beta):\alpha,
\beta \in V \cup U$ and $\alpha R^G \beta \vee \beta R^G \alpha\})$.  
The cardinality of $G$ is
$(|U^G|,|V^G|)$ or $|U^G| + |V^G|$. 

\noindent
2) We say $f$ is a strong embedding of the bipartite graph $G_1$ into the
bipartite graph $G_2$ \If \,:
\mn
\begin{enumerate}
\item[$(a)$]   $f$ is a one-to-one function from $U^{G_1} \cup V^{G_1}$
into $U^{G_2} \cup V^{G_2}$ mapping $U^{G_1}$ into $U^{G_2}$ and mapping
$V^{G_1}$ into $V^{G_2}$
\sn
\item[$(b)$]   for $(\alpha,\beta) \in U^{G_1} \times V^{G_1}$ we have

$\alpha R^{G_1} \beta \Leftrightarrow f(\alpha) R^{G_2} f(\beta)$.  
\end{enumerate}
\mn
3)  We say $f$ is a weak embedding of the bipartite graph $G_1$ into the
bipartite graph $G_2$ if
\mn
\begin{enumerate}
\item[$(a)$]   $f$ is a one-to-one function from $U^{G_1} \cup V^{G_1}$
into $U^{G_2} \cup V^{G_2},f$ mapping $U^{G_1}$ into $U^{G_2}$ and mapping
$V^{G_1}$ into $V^{G_2}$
\sn
\item[$(b)$]   for $(\alpha,\beta) \in U^{G_1} \times V^{G_1}$ we have

$\alpha R^{G_1} \beta \Rightarrow f(\alpha) R^{G_2} f(\beta)$.  
\end{enumerate}
4) In parts (2), (3) above, if $G_1$ is a bipartite graph and $G_2$ is a
graph then we mean $G^{[\text{gr}]}_1,G_2$. 

\noindent
5)  The $(\kappa,\theta)$-complete bipartite graph $K_{\kappa,\theta}$ is
$(U,V,R)$ with $U = \{i:i < \kappa\},V = \{\kappa +i:i < \theta\},
R = \{(i,\kappa +j):i < \kappa,j < \theta\}$, or any graph isomorphic
to it.
\end{definition}

\begin{definition}
\label{0.3}  
1) For a family ${\gH}$ of graphs
(or of bipartite graphs) we say $G$ is ste-universal [or we-universal]
for ${\gH}$ \Iff \,  every $G' \in {\gH}$ can be strongly embedded
[or weakly embedded] into $G$. 

\noindent
2) We say ${\gH}$ has a ste-universal (or we-universal) \If  \,  some
$G \in {\gH}$ is ste-universal (or we-universal) for ${\gH}$.
\end{definition}

\begin{definition}
\label{0.4}  
1) Let ${\gH}_{\lambda,\theta,\kappa} = {\gH}^{\text{gr}}_{\lambda,\theta,
\kappa}$ be the family of graphs $G$ of cardinality $\lambda$ (i.e. with
$\lambda$ nodes) such that the complete 
$(\theta,\kappa)$-bipartite graph cannot be weakly embedded into it;
gr stands for graph. 

\noindent
2) Let ${\gH}^{\text{bp}}_{{\bar \lambda},
\theta,\kappa}$ be the family of
bipartite graphs $G$ of cardinality $\bar \lambda$ such that the complete
$(\theta,\kappa)$-bipartite graph cannot be weakly embedded into it.
If $\bar \lambda = (\lambda,\lambda)$ we may write $\lambda$ (similarly in
(3)); bp stands for bipartite. 

\noindent
3) Let ${\gH}^{\text{sbp}}_{\bar \lambda,\{\theta,\kappa\}} = 
{\gH}^{\text{sbp}}_{\bar \lambda,\theta,\kappa}$ be the family
of bipartite graphs $G$ of cardinality $\bar \lambda$ such that
$K_{\theta,\kappa}$ (the $(\theta,\kappa)$-complete bipartite graph) and
$K_{\kappa,\theta}$ (the $(\kappa,\theta)$-complete bipartite graph) cannot
be weakly embedded into it; sbp stands for symmetrically bipartite.

\noindent
4) ${\gH}_\lambda = {\gH}^{\text{gr}}_\lambda$ is 
the family of graphs of cardinality $\lambda$ and 
${\gH}^{\text{bp}}_{\bar \lambda}$ is 
the family of bipartite graphs of cardinality $\bar \lambda$. 
\end{definition}

\begin{observation}
\label{0.5}  
1) The following are equivalent:
\mn
\begin{enumerate}
\item[$(a)$]   in ${\gH}^{\text{\rm sbp}}_{\lambda,\theta,\kappa}$
there is a we-universal
\sn
\item[$(b)$]   in $\{G^{[\text{\rm gr}]}:
G \in {\gH}^{\text{\rm sbp}}_{\lambda,\theta,\kappa}\}$
there is a we-universal.
\end{enumerate}
\mn
2) Similarly for ste-universal. 

\noindent
3) If $G$ is ste-universal for ${\frak H}$ then it is we-universal for
${\gH}$ (in all versions). 

\noindent
4) Assume that for every $G \in {\gH}_{\lambda,\theta,\kappa}$
there is a bipartite graph from ${\gH}_{\lambda,\theta,\kappa}$ not
$x$-embeddable into it, \then \,  in 
${\gH}^{\text{\rm bp}}_{\lambda,\kappa,\theta}$ and in
${\gH}^{\text{\rm bp}}_{\lambda,\theta,\kappa}$ and 
in ${\gH}^{\text{\rm sbp}}_{\lambda,\theta,\kappa}$
there is no x-universal member; for $x \in \{\text{we,\,ste}\}$. 
\end{observation}{

\begin{PROOF}{\ref{0.5}}
 (1) \underline{$(a) \Rightarrow (b)$}:  Trivially.
\medskip

\noindent
\underline{$(b) \Rightarrow (a)$}:  Assume $G$ is we-universal in
$\{G^{[\text{gr}]}:G \in {\gH}^{\text{sbp}}_{\lambda,\theta,\kappa}\}$ and 
let $\langle A_i:i < i^* \rangle$ be its connectivity
components.  Let $A_i$ be the disjoint union of $A_{i,0},A_{i,1}$ with no
$G$-edge inside $A_{i,0}$ and no $G$-edge inside $A_{i,1}$ (exists as $G =
G^{[\text{gr}]}_*$ for some $G_* \in 
{\frak H}^{\text{sbp}}_{\lambda,\kappa,\theta}$, note 
that $\{A_{i,0},A_{i,1}\}$ is unique
as $G \restriction A_i$ is connected).
Let $A^{m,\alpha}_{i,\ell}$ for $i < i^*,\ell < 2,m < 2,\alpha <
\lambda$ be pairwise disjoint sets with $|A^{m,\alpha}_{i,\ell}| =
|A_{i,k}|$ when $m + \ell = k$ mod $2$.  Let $G'$ be the following
member of ${\frak H}^{\text{bp}}_{\lambda,\theta,\kappa}$: let 
$U^{G'}$ be the disjoint union  of $A^{\ell,\alpha}_{i,\ell}$ for $i <
i^*,\ell < 2,\alpha < \lambda$ and $V^{G'}$ be the disjoint union of
$A^{1-\ell,\alpha}_{i,\ell}$ for $i < i^*,\ell < 2,\alpha < \lambda$
and $R^{G'} = \cup\{R^\alpha_{i,\ell}:i < i^*,\ell < 2\}$ where
$R^\alpha_{i,\ell}$ are chosen such that
$(A^{0,\alpha}_{i,0},
A^{1,\alpha}_{i,1},R^\alpha_{i,0}) \cong (A_{i,0},A_{i,1},R^G \restriction
A_{i,0} \times A_{i,1}) \cong
(A^{0,\alpha}_{i,1},A^{1,\alpha}_{i,0},R^\alpha_{i,1})$.  Easily $G'
\in {\frak H}^{\text{bp}}_{\lambda,\theta,\kappa}$ is we-universal. 

\noindent
2) The same proof.

\noindent
3), 4)  Easy. 
\end{PROOF}

\section {Some no we-universal} 

We show that if $\lambda = \lambda^\kappa \wedge 2^\kappa \ge \theta
\ge \kappa$ then in ${\gH}_{\lambda,\theta,\kappa}$ there is no
ste-universal graph (in \ref{1.4}); for we-universal there is a
similar theorem if $\theta = \kappa^+$, Pr$(\lambda,\kappa)$, see
\ref{1.3} + Definition \ref{1.2}, (this holds when $\lambda = \lambda^\kappa =
\text{ cf}(\lambda),\clubsuit_{S^\lambda_\kappa}$).

\begin{claim}
\label{1.1}   Assume $\kappa$ is regular and 
$\clubsuit_{S^{\kappa^+}_\kappa}$ (see Definition \ref{1.2} 
below).  \Then \, there is no we-universal in ${\gH}_{\kappa^+,
\kappa^+,\kappa}$.
\end{claim}

\begin{definition}
\label{1.2}  
1) For regular $\kappa < \lambda$ let
$S^\lambda_\kappa = \{\delta < \lambda:
\text{cf}(\delta) = \text{ cf}(\kappa)\}$. 

\noindent
2) For regular $\lambda$ and stationary subset $S$ of $\lambda$ let
$\clubsuit_S$ means that for some $\bar A = \langle A_\delta:\delta \in S,
\delta \text{ limit}\rangle$ we have
\mn
\begin{enumerate}
\item[$(a)$]   $A_\delta$ is an unbounded subset of $\delta$
\sn
\item[$(b)$]   if $A$ is an unbounded subset of $\lambda$ then for some
(equivalently stationarily many) $\delta \in S$ we have $A_\delta \subseteq
A$.
\end{enumerate}
\mn
2A) For $\kappa < \lambda$ let $\stick_{\lambda,\kappa}$ mean 
\underline{that} for
some family ${\cA} \subseteq [\lambda]^\kappa$ of cardinality
$\lambda$ we have $(\forall B \in [\lambda]^\lambda)(\exists A \in
{\cA})(A \subseteq B)$.

\noindent
3) Pr$(\lambda,\kappa)$ for cardinals $\lambda > \kappa$ \underline{means} that
some ${\cF}$ exemplifies it, which means
\mn
\begin{enumerate}
\item[$(a)$]  ${\cF}$ is a family of $\le \lambda$ functions
\sn
\item[$(b)$]   every $f \in {\cF}$ is a partial function from $\lambda$
to $\lambda$
\sn
\item[$(c)$]    if $f \in {\cF}$ then $\kappa = 
\text{ otp}(\text{Dom}(f))$ and $f$ strictly increasing
\sn
\item[$(d)$]   $f \ne g \in {\cF} \Rightarrow \kappa >
|\text{Dom}(f) \cap \text{ Dom}(g)|$
\sn
\item[$(e)$]   if $g$ is a partial (strictly) increasing function 
from $\lambda$ to
$\lambda$ such that Dom$(g)$ has cardinality $\lambda$, \then \, $g$
extends some $f \in {\cF}$.
\end{enumerate}
\mn
4) Pr$'(\lambda,\delta)$ is defined similarly for $\delta$ a limit
ordinal but clauses (c) + (d) are replaced by
\mn
\begin{enumerate}
\item[$(c)'$]   if $f \in {\cF}$ \then \, 
$\delta = \text{ otp}(\text{Dom}(f))$ and $f$ is one to one
\sn
\item[$(d)'$]   if $f \ne g \in {\cF}$ then $\text{Dom}(f) \cap
\text{ Dom}(g)$ is a bounded subset of Dom$(f)$ and of Dom$(g)$.
\end{enumerate}
\end{definition}

\begin{observation}
\label{1.2A} 
1) We have
\mn
\begin{enumerate}
\item[$(i)$]    $\kappa = \text{\rm cf}(\kappa)$ and 
$\clubsuit_{S^\lambda_\kappa} \Rightarrow \text{\rm Pr}(\lambda,\kappa)$
\sn
\item[$(ii)$]   Pr$'(\lambda,\kappa) \Rightarrow
\text{ Pr}(\lambda,\kappa) \Rightarrow \stick_{\lambda,\kappa}$
\sn
\item[$(iii)$]    for any cardinal $\kappa$ we have Pr$(\lambda,\kappa)
\Leftrightarrow \text{ Pr}'(\lambda,\kappa)$.
\end{enumerate}
\mn
2)  If we weaken clause (c) of \ref{1.2}(3) to
\mn
\begin{enumerate}
\item[$(c)^-$]   $f \in {\cF} \Rightarrow |\text{Rang}(f)| =
\kappa = |\text{Dom}(f)|$ 
\end{enumerate}
\mn
we get equivalent statement (can combine with \ref{1.2A}(3)).

\noindent
3) The ``one to one'' in Definition \ref{1.2}(4), clause $(c)'$ 
is not a serious demand, that is, omitting it we get an equivalent 
definition.
\end{observation}

\begin{PROOF}{\ref{1.2A}}
Easy.

\noindent
1) E.g., clause $(iii)$ holds because for any one to
one $f:\kappa \rightarrow \text{ Ord}$, for some $A \in
[\kappa]^\kappa$ the function $f \restriction A$ is strictly
increasing (note it first to regular $\kappa$). 

\noindent
2) Left to the reader.

\noindent
3) Why?   Let pr:$\lambda \times \lambda \rightarrow \lambda$ be 1-to-1
onto and pr$_1$,pr$_2:\lambda \rightarrow \lambda$ be such that
$\alpha = \text{ pr}(\text{pr}_1(\alpha),\text{pr}_2(\alpha))$.

\noindent
Let ${\cF}$ be as in the Definition \ref{1.2}(4) old version.
Then $\{\text{pr}_1 \circ f:f \in {\cF}\}$ will exemplify the new
version, i.e. without the 1-to-1.

For the other direction, just take $\{f \in {\cF}:f \text{ is one
to one}\}$. 
\end{PROOF}

\begin{PROOF}{\ref{1.1}}
\underline{Proof of \ref{1.1}}  It follows from \ref{1.3} proved below as
$\clubsuit_{S^{\kappa^+}_\kappa}$ easily implies Pr$(\kappa^+,\kappa)$
for regular $\kappa$ (see \ref{1.2A}(1)).  
\end{PROOF}

\begin{claim}
\label{1.3}  
If {\rm Pr}$(\lambda,\kappa)$, so $\lambda > \kappa$
\then \, in ${\gH}_{\lambda,\kappa^+,\kappa}$ there is no we-universal.

\noindent
Moreover, for every $G^* \in {\gH}_{\lambda,\kappa^+,\kappa}$ there is
a bipartite $G \in {\gH}_{\lambda,\kappa^+,\kappa}$ 
of cardinality $\lambda$ not we-embeddable into it.
\end{claim}

\begin{PROOF}{\ref{1.3}}
Let $G^*$ be a given graph from ${\gH}_{\lambda,\kappa^+,
\kappa}$; \wilog \, $V^{G^*} = \lambda$.  

For any $A \subseteq \lambda$ let
\mn
\begin{enumerate}
\item[$(*)_0$]   $(a) \quad Y^0_A 
=: \{\beta < \lambda:\beta$ is $G^*$-connected
with every $\gamma \in A\}$
\sn
\item[${{}}$]   $(b) \quad 
Y^2_A =: \{\beta < \lambda:\beta$ is $G^*$-connected with
$\kappa$ members of $A\}$
\sn
\item[${{}}$]   $(c) \quad Y^1_A 
=: \{\beta < \lambda:\beta$ is $G^*$-connected with
every $\gamma \in A$ except possibly 

\hskip25pt $< \kappa$ of them$\}$. 

Clearly
\sn 
\item[${{}}$]   $(d) \quad A \subseteq B \subseteq \lambda 
\Rightarrow Y^0_A \supseteq Y^0_B$ and
$Y^2_A \subseteq Y^2_B$ and $Y^1_A \supseteq Y^1_B$
\sn
\item[${{}}$]    $(e) \quad |A| \ge \kappa \Rightarrow Y^0_A 
\subseteq Y^1_A \subseteq Y^2_A$.
\end{enumerate}
\mn
We now note
\mn
\begin{enumerate}
\item[$(*)_1$]  if $A \in [\lambda]^\kappa$ then $|Y^0_A| \le \kappa$.
\end{enumerate}
\mn
[Why?  Otherwise we can find a weak embedding of the 
$(\kappa,\kappa^+)$-complete bipartite graph into $G^*$]
\mn
\begin{enumerate}
\item[$(*)_2$]  if $A \in [\lambda]^\kappa$ then $|Y^1_A| \le \kappa$.
\end{enumerate}
\mn
[Why?  If not choose pairwise disjoint subsets $A_i$ of $A$ for $i <
\kappa$ each of cardinality $\kappa$, now easily $\gamma \in Y^1_A 
\Rightarrow |\{i < \kappa:\gamma \notin Y^0_{A_i}\}| < \kappa$ 
so $Y^1_A \subseteq \bigcup\limits_{i < \kappa} Y^0_{A_i}$ hence
if $|Y^1_A| > \kappa$ then for some $i < \kappa,Y^0_{A_i}$ has
cardinality $> \kappa$, contradiction by $(*)_1$.]

Let ${\cF} = \{f_\alpha:\alpha < \lambda\}$ exemplify
Pr$(\lambda,\kappa)$.  
Now we start to choose the bipartite graph $G$:
\mn
\begin{enumerate}
\item[$\boxtimes_0$]   $U^G = \lambda,V^G = \lambda \times \lambda,R^G =
\bigcup\limits_{\alpha < \lambda} R^G_\alpha$ and $R^G_\alpha \subseteq
\{(\beta,(\alpha,\gamma)):\alpha < \lambda$ and
$\beta \in \text{ Dom}(f_\alpha)$ and $\gamma <
\lambda\} \subseteq U^G \times V^G$ where $R^G_\alpha$ is chosen below; we
let $G_\alpha = (U^G,V^G,R^G_\alpha)$.
\end{enumerate}
\mn
Now
\mn
\begin{enumerate}
\item[$\boxtimes_1$]   $G$ is a bipartite graph of cardinality $\lambda$
\sn
\item[$\boxtimes_2$]   the $(\kappa^+,\kappa)$-complete bipartite
graph $(\in {\gH}^{\text{bp}}_{(\kappa^+,\kappa)})$
cannot be weakly embedded into $G$.
\end{enumerate}
\mn
[Why?  As for any $(\alpha,\gamma) \in V^G$ the set $\{\beta <
\lambda:\beta R^G(\alpha,\gamma)\}$ is equal to Dom$(f_\alpha)$
which has cardinality $\kappa$ which is $< \kappa^+$;
 note that we are speaking of weak embedding as
bipartite graph, ``side preserving"]
\mn
\begin{enumerate}
\item[$\boxtimes_3$]   the $(\kappa,\kappa^+)$-complete bipartite
graph cannot be weakly embedded into $G$ provided that for each 
$\alpha < \lambda$,
\begin{enumerate}
\item[$\oplus^3_\alpha$]   $K_{\kappa,\kappa^+}$
cannot be weakly embedded into $(U^G,V^G,R^G_\alpha)$.
\end{enumerate}
\end{enumerate}
\mn
[Why?  Let $U_1 \subseteq U^G,V_1 \subseteq V^G$ have cardinality
$\kappa,\kappa^+$ respectively and let $V'_1 =
\{\alpha:(\alpha,\gamma) \in V_1$ for some $\gamma\}$.  If $|V'_1| \ge
2$ choose $9\alpha_1,\gamma_1),(\alpha_2,\gamma_2) \in V_1$ such that
$\alpha_1 \ne \alpha_2$, so $\{\beta <
\lambda:(\beta,(\alpha_\ell,\gamma_\ell)) \in R^G$ for $\ell=1,2\}$
include $U_1$ hence by the definition of $R^G$ we have
$|\text{Dom}(f_{\alpha_1}) \cap \text{ Dom}(f_{\alpha_2})| \ge |U_1| =
\kappa$, but 
$\alpha_1 \ne \alpha_2 \Rightarrow |\text{Dom}(f_{\alpha_1}) \cap
\text{Dom}(f_{\alpha_2})| < \kappa$ by clause (d) of Definition
\ref{1.2}(3) - think.]
\mn
\begin{enumerate}
\item[$\boxtimes_4$]  $G$ cannot be weakly embedded into $G^*$
provided that for each $\alpha < \lambda$:
\begin{enumerate}
\item[$\oplus^4_\alpha$]   there is no weak embedding 
$f$ of $(U^G,V^G,R^G_\alpha)$ into $G^*$ extending $f_\alpha$.
\end{enumerate}
\end{enumerate}
\mn
[Why?  Assume toward contradiction that $f$ is a one-to-one mapping
from $U^G \cup V^G$ into $V^{G^*} = \lambda$ mapping edges of $G$ to
edges of $G^*$.  So $f \restriction U^G$ is a one-to-one mapping from
$\lambda$ to $\lambda$ hence 
by the choice of ${\cF} = \{f_\alpha:\alpha < \lambda\}$ to witness
Pr$(\lambda,\kappa)$ see clause (e) of Definition \ref{1.2}(3)
there is $\alpha$ such that $f_\alpha \subseteq f \restriction U^G$.
So clearly $\beta R^G(\alpha,\gamma)$ and $\beta \in \text{
Dom}(f_\alpha)$ implies $\{f(\beta),f(\alpha,\gamma)\} =
\{f_\alpha(\beta),f((\alpha,\gamma))\} \in E^{G^*}$, hence
$(\beta,(\alpha,\gamma)) \in R^G_\alpha \Rightarrow
\{f(\beta),f((\alpha,\gamma))\} \in E^{G^*}$.  This clearly
contradicts $\otimes^4_\alpha$ which we are assuming.]

So we are left with, for each $\alpha < \lambda$, choosing $R_\alpha
\subseteq \{(\beta,(\alpha,\gamma):\beta \in \text{ Dom}(f_\alpha),\gamma <
\lambda\}$ to satisfy $\oplus^3_\alpha + \oplus^4_\alpha$.  
The proof splits to cases, fixing $\alpha$.

Let us denote $B_\alpha = \text{ Rang}(f_\alpha),A_\alpha =
\text{ Dom}(f_\alpha)$ for $\ell = 0,1$ we let 
$A^\ell_\alpha =: \{\gamma \in A_\alpha:\text{otp}
(A_\alpha \cap \gamma) = \ell \text{ mod } 2\}$ and $B^\ell_\alpha
=: \{f_\alpha(\gamma):\gamma \in A^\ell_\alpha\}$.

\noindent
\underline{Case 1}:  $Y^2_B$ has cardinality $\le \kappa$ for some $B \in
[B_\alpha]^\kappa$.

Choose such $B = B'_\alpha$ and let $A'_\alpha = \{\beta \in A_\alpha:
f_\alpha(\beta) \in B'_\alpha\}$.  There is a sequence $\bar C = \langle
C_\zeta:\zeta < \kappa^+ \rangle,C_\zeta \in [\kappa]^\kappa$ such that
$\xi < \zeta \Rightarrow |C_\xi \cap C_\zeta| < \kappa$.

Let $R_\alpha = \{(\beta,(\alpha,\gamma)):\gamma < \kappa^+,\beta \in
A'_\alpha$ and otp$(\beta \cap A'_\alpha) \in C_\gamma\}$.  Now
$\oplus^4_\alpha$ holds because if $f$ is a counter-example, \then \,
necessarily by the pigeon-hole principle 
for some $\gamma < \kappa^+$ we have
$f((\alpha,\gamma)) \notin Y^2_{B'_\alpha}$, but clearly $(\alpha,\gamma)$ is
$G_\alpha$-connected to $\kappa$ members of $A'_\alpha$ hence $f((\alpha,
\gamma))$ is $G^*$-connected to $\kappa$ members of $B'_\alpha$ hence
$f((\alpha,\gamma)) \in Y^2_{B'_\alpha}$ and we get a contradiction.
Also $\oplus^3_\alpha$ holds as $\xi < \zeta < \kappa^+ \Rightarrow
|C_\xi \cap C_\zeta| < \kappa$.

So we may assume, for the rest of the proof, that
\mn
\begin{enumerate}
\item[$\boxtimes_5$]   $|Y^2_B| > \kappa$ for every $B \in
[B_\alpha]^\kappa$.
\end{enumerate}
\bigskip

\noindent
\underline{Case 2}:  For some 
$\ell < 2,|Y^2_{B^\ell_\alpha}| > \kappa$ and for some $Z$:
\mn
\begin{enumerate}
\item[$(i)$]   $Z \subseteq Y^2_{B^\ell_\alpha} \backslash
Y^1_{B^\ell_\alpha}$
\sn
\item[$(ii)$]   $|Z| \le \kappa$
\sn
\item[$(iii)$]  for every $\gamma_0 \in (Y^2_{B^\ell_\alpha} \backslash
Y^1_{B^\ell_\alpha}) \backslash Z$ there is $\gamma_1 \in Z$ such that 
$\kappa > |\{\beta \in B^\ell_\alpha:\beta \text{ is } G^*
\text{-connected to } \gamma_0$ but is not $G^*$-connected to $\gamma_1\}|$.
\end{enumerate}
\mn
So we choose such $\ell = \ell(\alpha) < 2,Z = Z_\alpha$ and then
we choose a sequence $\langle B_{\alpha,\gamma}:\gamma \in
Z_\alpha \rangle$ such that:
\mn
\begin{enumerate}
\item[$\boxtimes_6$]  $(a) \quad B_{\alpha,\gamma}$ 
is a subset of $B^\ell_\alpha$
\sn
\item[${{}}$]  $(b) \quad |B_{\alpha,\gamma}| = \kappa$
\sn
\item[${{}}$]   $(c) \quad \gamma_1 \ne \gamma_2 \in Z \Rightarrow 
B_{\alpha,\gamma_1} \cap B_{\alpha,\gamma_2} = \emptyset$
\sn
\item[${{}}$]  $(d) \quad \gamma$ is not $G^*$-connected to 
any $\varepsilon \in B_{\alpha,\gamma}$
\end{enumerate}
\mn
(this is possible as $\gamma \in Z_\alpha \Rightarrow \gamma
\notin Y^1_{B^\ell_\alpha}$).

Now we can find a sequence $\langle C_{\alpha,\zeta}:\zeta < \kappa^+
\rangle$ satisfying
\mn
\begin{enumerate}
\item[$\boxtimes_7$]  $(\alpha) \quad C_{\alpha,\zeta} 
\subseteq B^\ell_\alpha$
\sn
\item[${{}}$]  $(\beta) \quad |C_{\alpha,\zeta}| = \kappa$ 
moreover $\beta \in Z_\alpha \Rightarrow |C_{\alpha,\zeta} \cap 
B_{\alpha,\beta}| = \kappa$
\sn
\item[${{}}$]  $(\gamma) \quad$ for $\xi < \zeta$ we have 
$|C_{\alpha,\xi} \cap C_{\alpha,\zeta}| < \kappa$
\end{enumerate}
\mn
(e.g. if $\kappa = \text{ cf}(\kappa) > \aleph_0$ by renaming $B^\ell_\alpha
= \kappa$, each $B_{\alpha,\varepsilon}$ is stationary, choose nonstationary
$C_{\alpha,\varepsilon} \subseteq \kappa$ inductively on $\varepsilon$; if
$\kappa > \text{ cf}(\kappa)$ reduce it to construction on regulars, if
$\kappa = \aleph_0$ like $\kappa = \text{ cf}(\kappa) > \aleph_0$). 

Lastly we choose 
$R_\alpha = \{(\beta,(\alpha,\gamma)):\beta \in A_\alpha,\gamma
< \kappa^+ \text{ and } f_\alpha(\beta) \in C_{\alpha,\gamma}\}$. 

Now $\oplus^3_\alpha$ is proved as in the first case, as for 
$\oplus^4_\alpha$, if $f$ is a counter-example \then \,
clearly for $\gamma < \kappa^+,f((\alpha,\gamma)) \in Y^2_{B_\alpha}$,
so as $|Y^2_{B_\alpha}| > \kappa$ by $\boxtimes_5$ and
$|Z_\alpha| \le \kappa$ and $|Y^1_{B^\ell_\alpha}| 
\le \kappa$ (by $(*)_2$) necessarily for some
$\zeta < \kappa^+,\gamma_0 =: f((\alpha,\zeta)) 
\in Y^2_{B^\ell_\alpha} \backslash
Y^1_{B^\ell_\alpha} \backslash Z_\alpha$. Let $\gamma_1 \in Z_\alpha$ be as
guaranteed in clause (iii) in the present case.
Now $\gamma_0$ is $G^*$-connected to every member of
$C_{\alpha,\zeta}$ as $\gamma_0 = f((\alpha,\zeta))$.  Hence
$\gamma_0$ is $G^*$-connected to $\kappa$ members of $B_{\alpha,\gamma_1}$
(see clause $(\beta)$ above and the choice of $R_\alpha$); 
but $\gamma_1$ is not $G^*$-connected to
any member of $B_{\alpha,\gamma_1}$ (see clause (d) above).  Reading
clause (iii), we get contradiction.
\bigskip

\noindent
\underline{Case 3}:  \underline{Neither Case 1 nor Case 2}.

Recall that $\alpha$ is fixed.
For $\ell \in \{0,1\}$ we choose $Z^\ell_{\alpha,\zeta}$ by induction on
$\zeta < \kappa^+$, such that
\mn
\begin{enumerate}
\item[$\boxtimes_8$]   $(a) \quad 
Z^\ell_{\alpha,\zeta}$ a subset of $Y^2_{B^\ell_\alpha}$ 
of cardinality $\kappa$
\sn
\item[${{}}$]  $(b) \quad Z^\ell_{\alpha,\zeta}$ is increasing continuous in
$\zeta$
\sn
\item[${{}}$]  $(c) \quad Y^1_{B^\ell_\alpha} \subseteq Z^\ell_{\alpha,0}$
\sn
\item[${{}}$]  $(d) \quad$ if $\zeta = \xi +1$ \then \, there is 
$\gamma^\ell_{\alpha,\xi} \in Z^\ell_{\alpha,\zeta} \backslash 
Z^\ell_{\alpha,\xi} \backslash Y^1_{B^\ell_\alpha}$ such that 

\hskip25pt for every $\gamma' \in Z^\ell_{\alpha,\xi} \backslash 
Y^1_{B^\ell_\alpha}$ we have 

\hskip25pt $\kappa = |\{\beta \in B^\ell_\alpha:
\beta \text{ is } G^*$-connected to
$\gamma^\ell_{\alpha,\xi}$ but not to $\gamma'\}|$
\sn
\item[${{}}$]  $(e) \quad$ if $\zeta = \xi +1$ and $\gamma \in
Z^\ell_{\alpha,\zeta}$ hence $\kappa = |\{\beta \in B^\ell_\alpha:\beta$
is connected to $\gamma\}|$ 

\hskip25pt (e.g. $\gamma =
\gamma^\ell_{\alpha,\xi}$), \then \,  
$Y^1_{\{\beta \in B^\ell_\alpha:\beta\text{ is } G^*\text{-connected
to }\gamma\}}$ is included in 

\hskip25pt $Z^\ell_{\alpha,\zeta}$
\sn
\item[${{}}$]  $(f) \quad Z^0_{\alpha,\zeta} \cap (Y^2_{B^0_\alpha} \cap 
Y^2_{B^1_\alpha}) = Z^1_{\alpha,\zeta} \cap (Y^2_{B^0_\alpha}
\cap Y^2_{B^1_\alpha})$.
\end{enumerate}
\mn
Why possible?  For clause (c) we have $|Y^1_{B^\ell_\alpha}| \le
\kappa$ by $(*)_2$, for clause (d) note that ``not Case 2"
trying $Z^\ell_{\alpha,\xi} \backslash Y^1_{B^\ell_\alpha}$ as $Z$, and for
clause (e) note again 
$|Y^1_{\{\beta \in B^\ell_\alpha:\beta \text{ is }
G^*\text{-connected to } \gamma^\ell_{\alpha,\varepsilon}\}}| \le
\kappa$ by $(*)_2$.

Having chosen $\langle Z^\ell_{\alpha,\zeta}:\zeta < \kappa^+,\ell < 2
\rangle$, we let

\begin{equation*}
\begin{array}{clcr}
R_\alpha = \{(\beta,(\alpha,\zeta)):&\text{ for some } \ell < 2
\text{ we have:} \\
  &\beta \in A^\ell_\alpha \text{ and } f_\alpha(\beta) \text{ is }
G^* \text{-connected to } \gamma = \gamma^\ell_{\alpha,2 \zeta + \ell}
\text{ so } \zeta < \kappa^+\}.
\end{array}
\end{equation*}
\mn
Now why $\oplus^3_\alpha$ holds?  Otherwise, we can find $A \subseteq
A_\alpha,|A| = \kappa$ and $B \subseteq \kappa^+,|B| = \kappa^+$ such that
$\beta \in A$ and $\xi \in B \Rightarrow 
\beta R_\alpha(\alpha,\gamma^{\ell(\xi)}_{\alpha,\xi})$ where $\xi =
\ell(\xi)$ mod $2$, so for some 
$\ell < 2$ we have $|A \cap A^\ell_\alpha| = \kappa$, and let

\[
A' = \{f_\alpha(\beta):\beta \in A \cap A^\ell_\alpha\}.
\]
\mn
Easily $|B| = \kappa^+$ and $|A'| = \kappa$ and 
$\beta \in A'$ and $\xi \in B \Rightarrow \beta R^{G^*}
\gamma^\ell_{\alpha,\xi}$, contradiction to 
``$K_{\kappa,\kappa^+}$ is not weakly embeddable into $G^*$".

Lastly, why $\oplus^4_\alpha$ holds?  Otherwise, letting $f$ be a
counterexample, let $\zeta < \kappa^+$ and $\ell < 2$.  Clearly
$f((\alpha,\zeta))$ 
is $G^*$-connected to every $\beta' \in B^\ell_\alpha$ which is
$G^*$-connected to $\gamma^\ell_{\alpha,2 \zeta + \ell}$ hence 
$f((\alpha,\zeta))$
cannot belong to $Z^\ell_{\alpha,2 \zeta + \ell}
\backslash Y^1_{B^\ell_\alpha}$ (by the demand in clause (d) of 
$\boxtimes_8$), but it has to belong to
$Z^\ell_{\alpha,2 \zeta + \ell +1}$ (by clause (e) of $\boxtimes_8$),
so $f((\alpha,\zeta)) \in (Z^\ell_{\alpha,2 \zeta + \ell +1} \backslash
Z^\ell_{\alpha,2 \zeta + \ell}) \cup Y^1_{B^\ell_\alpha}$.  Putting together
$\ell = 0,1$ we get $f((\alpha,\zeta)) \in ((Z^0_{\alpha,2 \zeta +1}
\backslash Z^0_{\alpha,2 \zeta}) \cup Y^1_{B^0_\alpha}) \cap
((Z^1_{\alpha,2 \zeta +2} \backslash Z^1_{\alpha,2 \zeta +1}) \cup
Y^1_{B^1_\alpha})$ hence $f((\alpha,\zeta)) \in
Y^1_{B^0_\alpha} \cup Y^1_{B^1_\alpha}$,
but $|Y^1_{B^\ell_\alpha}| < \kappa^+$; as this holds for every $\zeta
< \kappa^+$ this is a contradiction to ``$f$ is one to one''.
\end{PROOF}

\begin{claim}
\label{1.4}  
1)  Assume $\lambda \ge 2^\kappa \ge \theta \ge
\kappa$ and $\lambda = \lambda^\kappa$ (e.g. $\lambda = 2^\kappa$). 

\Then \, in ${\gH}_{\lambda,\kappa,\theta}$ there is 
no {\rm ste}-universal (moreover, the counterexamples are bipartite). 

\noindent
2) Assume {\rm Pr}$(\lambda,\kappa),\lambda \ge \theta \ge \kappa,2^\kappa
\ge \theta$.  \Then \, the conclusion of (1) holds.
\end{claim}

\begin{PROOF}{\ref{1.4}}
1) By the simple black box (\cite[Ch.III,\S4]{Sh:300}) or
\cite[Ch.VI,\S1]{Sh:e}, i.e. \cite{Sh:309})
\mn
\begin{enumerate}
\item[$\boxtimes$]   there is $\bar f = \langle f_\eta:\eta \in
{}^\kappa \lambda \rangle,f_\eta$ a function from $\{\eta \restriction
i:i < \kappa\}$ into $\lambda$ such that for every $f:{}^{\kappa >}
\lambda \rightarrow \lambda$ for some $\eta \in {}^\kappa \lambda$ we have
$f_\eta \subseteq f$.
\end{enumerate}
\mn
Let $G^* \in {\gH}_{\lambda,\kappa,\theta}$ and we shall show that it
is not ste-universal in ${\gH}_{\lambda,\kappa,\theta}$, without
loss of generality $V^{G^*} = \lambda$.  For this we
define the following bipartite graph $G$:
\mn
\begin{enumerate}
\item[$\boxplus_1$]   $(i) \quad U^G = 
{}^{\kappa >}\lambda$ and $V^G = {}^\kappa \lambda$
\sn
\item[${{}}$]  $(ii) \quad R^G = \cup\{R^G_\eta:
\eta \in {}^\kappa \lambda$ and $f_\eta$ is a one-to-one function$\}$ where

\hskip25pt $R^G_\eta \subseteq \{(\eta \restriction i,\eta):i \in u_\eta\}$ 
where $u_\eta \subseteq \kappa$ is defined as follows
\sn
\item[$\boxplus_2$]  for $\eta \in {}^\kappa \lambda$
we choose $u_\eta \subseteq \kappa$ such that if possible
\begin{enumerate}
\item[${{}}$]  $(*)_{\eta,u_\eta} \quad$ for no $\gamma < \lambda$ do
we have $(\forall i < \kappa)[f_\eta(\eta \restriction i) R^{G^*} \gamma
\equiv i \in u_\eta]$.
\end{enumerate}
\end{enumerate}
\mn
If for every $\eta \in {}^\kappa \lambda$ for which $f_\eta$ is one to
one for some $u \subseteq \kappa$ we have 
$(*)_{\eta,u}$ holds, \underline{then} clearly by $\boxtimes$ we are done.

Otherwise, for this $\eta \in {}^\kappa \lambda,f_\eta$ is one to one
and: there is $\gamma_u <
\lambda$ satisfying $(\forall i < \kappa)(f_\eta(\eta \restriction i) R^{G^*}
\gamma_u \Leftrightarrow i \in u)$ for every $u \subseteq \kappa$.  But then
$A' =: \{f_\eta(\eta \restriction 2i):i < \kappa\}$ and $B' =: \{\gamma_u:
u \subseteq \kappa \text{ and } (\forall i < \kappa)2i \in u\}$ form a
complete $(\kappa,2^\kappa)$-bipartite subgraph of $G^*$,
contradiction. 

\noindent
2) The same proof. 
\end{PROOF}

\section {No we-universal by Pr$(\lambda,\kappa)$ and its relatives}

We define here some relatives of Pr.  Here Ps is like Pr but we are
approximating $f:\lambda \rightarrow \lambda$, and
Pr$_3(\chi,\lambda,\mu,\alpha)$ is a weak version of $(\lambda +
\mu)^{|\alpha|} \le \chi$ (Definition \ref{2.4A}); we give
sufficient conditions by cardinal arithmetic (Claim \ref{2.5},
\ref{2.5A}).  We prove more cases of no we-universal: the case
$\theta$ limit (and {\rm Pr}$(\lambda,\kappa)$) in \ref{2.2}, a case of
Pr$'(\lambda,\theta^+ \times \kappa)$ in \ref{2.4}.  We also note
that we can replace Pr by Ps in \ref{2.5B}, and $\lambda$ strong
limit singular of cofinality $> \text{ cf}(\kappa)$ in \ref{2.5C}.

\begin{convention}
\label{2.1}  $\lambda \ge \theta \ge \kappa \ge \aleph_0$.
\end{convention}

\begin{claim}
\label{2.2}   If $\theta$ is a limit cardinal and 
{\rm Pr}$(\lambda,\kappa)$, \then \, there is no {\rm we}-universal 
graph in ${\gH}_{\lambda,\theta,\kappa}$ even for the class 
of bipartite members.
\end{claim}

\begin{PROOF}{\ref{2.2}}   
Like the proof of \ref{1.3}, except that we 
replace cases 1-3 by:

for every $\alpha < \lambda$ 
we let $R_\alpha = \{(\beta,(\alpha,\gamma)):\beta \in \text{ Dom}(f_\alpha)$
and $\gamma < |Y^0_{\text{Dom}(f_\alpha)}|^+\}$.  Now $\oplus^3_\alpha$
holds as $|Y^0_{\text{Dom}(f_\alpha)}| < \theta$ (by $(*)_1$ there) hence
$|Y^0_{\text{Dom}(f_\alpha)}|^+ < \theta$ as $\theta$ is a limit cardinal.
Lastly $\oplus^4_\alpha$ holds as for some $\alpha$ we have $f_\alpha
\subseteq f$ hence the function $f$ maps $\{(\alpha,\gamma):\gamma <
|Y^0_{\text{Dom}(f_\alpha)}|^+\}$ into $Y^0_{\text{Dom}(f_\alpha)}$ but $f$
is a one to one mapping, contradiction.  
\end{PROOF}

Recall
\begin{definition}
\label{pr.22}
For a cardinal $\lambda$ and a limit ordinal $\delta$,
Pr$'(\lambda,\delta)$ holds \when \, for some $\cF$:
\mn
\begin{enumerate}
\item[$(a)$]  $\cF$ a family of $\le \lambda$ functions
\sn
\item[$(b)$]  every $f \in \cF$ is a partial function from $\lambda$
to $\lambda$
\sn
\item[$(c)$]  $f \in \cF \Rightarrow \text{ otp}(\text{Dom }f) =
\delta$, and $f$ is one to one
\sn
\item[$(d)$]  $f,g \in \cF,f \ne g \Rightarrow (\text{Dom }f) \cap
(\text{Dom }g)$ is a bounded subset of Dom$(f)$ and of Dom$(g)$
\sn
\item[$(e)$]  if $g:\lambda \rightarrow \lambda$ is a partial
function, one to one, and $|\text{Dom }g| = \lambda$, \then \, $g$
extends some $f \in \cF$.
\end{enumerate}
\end{definition}

\begin{claim}
\label{2.4}  Assume
\mn
\begin{enumerate}
\item[$(a)$]  {\rm Pr}$'(\lambda,\delta^*),
\delta^* = \sigma \times \kappa$,
ordinal product \footnote{this is preserved by decreasing $\sigma$}
\sn
\item[$(b)$]   $\sigma = \theta^+$.
\end{enumerate}
\mn
\Then \, there is no {\rm we}-universal in 
${\gH}^{\text{\rm gr}}_{\lambda,\theta,
\kappa}$ even for the class of bipartite members.
\end{claim}

\begin{PROOF}{\ref{2.4}}  
Let $G^* \in {\gH}_{\lambda,\theta,\kappa}$ and we shall
prove it is not we-universal; let without loss of generality 
$V^{G^*} = \lambda$.

Let ${\cF}$ be a family 
exemplifying Pr$'(\lambda,\delta^*)$, let ${\cF} =
\{f_\alpha:\alpha < \lambda\}$ let $A_\alpha = \text{ Dom}(f_\alpha)$ and
let it be $\{\beta_{\alpha,\varepsilon,i}:i < \sigma,\varepsilon < \kappa\}$
such that 
$[\beta_{\alpha,\varepsilon(1),i(1)} <
\beta_{\alpha,\varepsilon(2),i(2)} \Leftrightarrow 
\varepsilon(1) < \varepsilon(2) \vee (\varepsilon(1) =
\varepsilon(2)$ and $i(1) < i(2))]$ and for $i < \sigma$
let $A_{\alpha,i} = \{\beta_{\alpha,
\varepsilon,i}:\varepsilon < \kappa\}$, so clearly
\mn
\begin{enumerate}
\item[$(*)_1$]   $A_{\alpha,i} \in [\lambda]^\kappa$ and $(\alpha_1,i_1) \ne
(\alpha_2,i_2) \Rightarrow |A_{\alpha_1,i_1} \cap A_{\alpha_2,i_2}| <
\kappa$.
\end{enumerate}
\mn
[Why the second assertion?  As $\{\beta_{\alpha,\varepsilon,i}:\varepsilon <
\kappa\}$ is an unbounded subset of $A_\alpha$ (of order type $\kappa$).]

For $(\alpha,i) \in \lambda \times \sigma$ let $f_{\alpha,i} = f_\alpha
\restriction A_{\alpha,i}$ let $B_{\alpha,i} = \text{ Rang}(f_{\alpha,i})$
so $|A_{\alpha,i}| = |B_{\alpha,i}| = \kappa$ and let $Y^0_{\alpha,i} = 
\{\gamma < \lambda:\gamma$ is $G^*$-connected to every member of 
$B_{\alpha,i}\}$, so as $G^* \in {\frak H}_{\lambda,\theta,\kappa}$
clearly $|Y^0_{\alpha,i}| < \theta$.  As $\sigma = \theta^+ > \theta
\ge \kappa$, clearly for each
$\alpha < \lambda$ for some $\mu_\alpha < \theta$ we have \footnote{in
fact by \ref{2.2} \wilog \, $\theta$ is a successor cardinal, so \wilog \,
$\mu^+_\alpha = \theta$}: 
$X_\alpha := \{i < \sigma:|Y^0_{\alpha,i}| \le \mu_\alpha\}$
has cardinality $\sigma$.  As $\mu_\alpha < \theta$ also $\chi_\alpha
:= \mu^+_\alpha$ is $< \theta^+ = \sigma$, so $\chi^+_\alpha \le
\sigma = |X_\alpha|$. 
We choose by induction on $\varepsilon < \chi^+_\alpha$ an ordinal 
$i^*_{\alpha,\varepsilon} \in X_\alpha$ such that:
\mn
\begin{enumerate}
\item[$(*)_2$]   $i^*_{\alpha,\varepsilon} \notin
\{i^*_{\alpha,\zeta}:\zeta < \varepsilon\}$.
\end{enumerate}
\mn
Recall that $\chi_\alpha$ is a successor, hence a regular cardinal.
So if $\varepsilon \in S^{\chi^+_\alpha}_{\chi_\alpha} = \{\varepsilon
< \chi^+_\alpha:\text{cf}(\varepsilon) = \chi_\alpha\}$ recalling
$\chi_\alpha = \mu^+_\alpha > |Y^0_{\alpha,i^*_{\alpha,\varepsilon}}|$
there is $\zeta < \varepsilon$ such that 
$(Y^0_{\alpha,i^*_{\alpha,\varepsilon}} \cap \cup
\{Y^0_{\alpha,\xi}:\xi < \varepsilon\} \subseteq \cup
\{Y^0_{\alpha,\xi}:\xi < \zeta\})$.  Let $g(\varepsilon)$ be the first
ordinal having this property, so $g$ is a well defined function with
domain $S^{\chi^+}_{\chi_\alpha}$.  Clearly, $g$ is a regressive funciton.

By Fodor's lemma for some
stationary $S_\alpha \subseteq S^{\chi^+_\alpha}_{\chi_\alpha}$, and for some
$B^*_\alpha$ of cardinality $\le \chi_\alpha$ we have
$\varepsilon \in S_\alpha \Rightarrow B^*_\alpha \supseteq
Y^0_{\alpha,i^*_{\alpha,\varepsilon}} \cap \cup\{Y^0_{\alpha,j}:
j < i^*_{\alpha,\varepsilon}\}$, in fact: $B^*_\alpha = \cup
\{Y^0_{\alpha,i^*_{\alpha,\varepsilon}}:\varepsilon < \varepsilon^*\}$
where $g \restriction S_\alpha$ is constantly $\varepsilon^*$ is O.K.,
we can decrease $B^*_\alpha$ but immaterial here  

Now
\mn
\begin{enumerate}
\item[$(*)_3$]   for $\xi \ne \zeta$ from $S_\alpha$ 
there is no $\beta < \lambda$ such that:

$\beta$ is $G^*$-connected to every 
$\gamma \in B_{\alpha,i^*_{\alpha,\xi}}$

$\beta$ is $G^*$-connected to every 
$\gamma \in B_{\alpha,i^*_{\alpha,\zeta}}$ 

$\beta$ is not in $B^*_\alpha$.
\end{enumerate}
\mn
Let $\langle \zeta(\alpha,j):j < \chi^+_\alpha \rangle$ list
$S_\alpha$ in an increasing order.  Let $G$ be the bipartite graph

$$
U^G = \lambda
$$

$$
V^G = \lambda \times \sigma
$$

$$
R^G = \bigl\{ (\beta,(\alpha,\gamma)):\alpha < \lambda \text{ and }
\gamma < \chi^+_\alpha \text{ and }
\beta \in A_{\alpha,i^*_{\alpha,\zeta(\alpha,2 \gamma)}} \cup
A_{\alpha,i^*_{\alpha,\zeta(\alpha,2 \gamma +1)}} \bigr\}.
$$
\mn
\begin{enumerate}
\item[$\boxtimes_1$]  $G$ is a bipartite graph of cardinality
$\lambda$
\sn
\item[$\boxtimes_2$]  the $(\theta,\kappa)$-complete bipartite graph
$(\in \gH^{\text{\rm bp}}_{(\theta,\kappa)}$) cannot be weakly
embedded into $G$.
\end{enumerate}
\mn
[Why $\boxtimes_2$ holds?  So let $(\alpha(1),\gamma(1)) \ne
(\alpha(2),\gamma(2))$ belongs to $V^G$, the set set
$\{\beta \in U^G:\beta$ connected to $(\alpha(1),\gamma(2))$ and to
$(\alpha(2),\gamma(2))\}$ is included in
$\bigcup\limits_{\iota(1),\iota(2) \in \{0,1\}}
(A_{\alpha(1),i^*_{\zeta(\alpha(1),2 \gamma(1) + \iota(1)}} \cap 
A_{\alpha(2),i^*_{\zeta(\alpha(2),2 \gamma(2) + \iota(2)}})$ which is
the union of four sets each of cardinal $< \kappa$ (by $(*)_1$) hence
has cardinality $< \kappa$.]
\mn
\begin{enumerate}
\item[$\boxtimes_3$]  the $(\kappa,\theta)$-complete bipartite graph
cannot be weakly embedded into $G$.
\end{enumerate}
\mn
[Why?  Toward contradiction assume $U_1 \subseteq U^G,V_1 \subseteq
V^G$ have cardinality $\kappa,\theta$ respectively and $\beta \in U_1
\wedge (\alpha,\gamma) \in V_1 \Rightarrow \beta R(\alpha,\gamma)$.

Let $(\alpha,\gamma) \in V_1$, clearly $\beta \in U_1 \Rightarrow
\beta R^G(\alpha,\gamma) \Rightarrow \beta \in
A_{\alpha,i^*_{\alpha,\zeta(\alpha,2 \gamma)}} \cup
A_{\alpha,i^*_{\alpha,\zeta(\alpha,2 \gamma +1)}}$.  So 
$U_1 \subseteq A_{i^*_{\alpha,\zeta(\alpha,2
\gamma)}} \cup A_{i^*_{\alpha,\zeta(\alpha,2 \gamma +1)}}$.  

Now if $(\alpha_1,\gamma_1),(\alpha_2,\gamma_2) \in V_1$ and $\alpha_1
\ne \alpha_2$ then $i,j < \sigma \Rightarrow |A_{\alpha_1,i} \cap
A_{\alpha_2,j}| < \kappa$ by the choice of $\cF$, so necessarily for
some $\alpha_* < \lambda$ we have $V_1 \subseteq \{\alpha_*\} \times
\sigma$.  But if $(\alpha,\gamma_1) \ne (\alpha,\gamma_2) \in V_1$
then $(\alpha,\gamma_1),(\alpha,\gamma_2)$ has no common neighbour,
contradiction. 
\mn
\begin{enumerate}
\item[$\boxtimes_4$]  there is no weak embedding $f$ of $G$ into
$G^*$.
\end{enumerate}
\mn
[Why?  Toward contradiction assume that $f$ is such a weak embedding.
By the choice of $\cF$ and $\langle f_\alpha:\alpha < \lambda\rangle$
we can choose $\alpha < \lambda$ such that $f_\alpha \subseteq f \rest
U^G$.  As
$f$ is a weak embedding $\beta < \lambda \wedge \gamma < \lambda
\wedge \beta R^G (\alpha,\gamma) \Rightarrow
f(\beta) R^{G^*}f(\alpha,\gamma)$, hence $\beta \in \bigcup\limits_{i} 
A_{\alpha,i} = \text{ Dom}(f_\alpha) \wedge \gamma < \lambda \wedge
\beta R^G(\alpha,\gamma) \Rightarrow
f_\alpha(\beta)R^{G^*}f(\alpha,\gamma)$.  Hence if $\gamma < \sigma$
then $f(\alpha,\gamma)$ is $G^*$-connected to every $\beta \in
B_{\alpha,\zeta(\alpha,2 \gamma)} \cup B_{\alpha,\zeta(\alpha,2
\gamma)} \cup \{\beta_{\alpha,\gamma}\}$ hence $f(\alpha,\gamma) \in
Y^0_{\alpha,\zeta(\alpha,2 \gamma)} \cap Y^0_{\alpha,\zeta(\alpha,2
\gamma +1)}$ which implies that $f(\alpha,\gamma) \in B^*_\alpha$.
So the function $f$ maps the set $\{(\alpha,\gamma):\gamma < \sigma\}$
into $B^*_\alpha$.   But $f$ is a one-to-one function and 
$B^*_\alpha$ has cardinality $< \sigma$, contradiction.    
\end{PROOF}

\begin{definition}
\label{2.4A}  
1) For $\kappa < \lambda$ and
$\delta < \lambda$ we define Ps$(\lambda,\kappa)$ and
Ps$'(\lambda,\delta)$ similarly to the definition of
Pr$(\lambda,\kappa)$, Pr$'(\lambda,\delta)$ in Definition
\ref{1.2}(3),(4) except that we replace clause (e) by
\mn
\begin{enumerate}
\item[$(e)^-$]   if $g$ is a one to one function from $\lambda$ to
$\lambda$, \then \, $g$ extends some $f \in {\cF}$
\end{enumerate}
\mn
(so the difference is that Dom$(g)$ is required to be $\lambda$). 

\noindent
2) Let Pr$_3(\chi,\lambda,\mu,\alpha)$ means that for some ${\cF}$:
\mn
\begin{enumerate}
\item[$(a)$]   ${\cF}$ a family of partial functions from $\mu$
to $\lambda$
\sn
\item[$(b)$]   $|{\cF}| \le \chi$
\sn
\item[$(c)$]   $f \in {\cF} \Rightarrow \text{ otp}(\text{Dom}(f)) = \alpha$
\sn
\item[$(d)$]   if $g \in {}^\mu \lambda \Rightarrow (\exists f \in
{\cF})(f \subseteq g)$.
\end{enumerate}
\end{definition}

\begin{claim}
\label{2.5}  
1) Assume $\lambda$ is strong limit, $\lambda > \kappa$ and 
{\rm cf}$(\lambda) > { \text{\rm cf\/}}(\kappa)$.  

\noindent
\Then \,
{\rm Pr}$'(\lambda,\delta^*)$ holds if $\delta^* < \lambda$ has cofinality
{\rm cf}$(\kappa)$. 

\noindent
2) If $\lambda = \mu^+ = 2^\mu$, 
{\rm cf}$(\delta^*) \ne { \text{\rm cf\/}}(\mu),
\delta^* < \lambda$ \then \, {\rm Pr}$'(\lambda,\delta^*)$ holds. 

\noindent
3) If $\delta < \lambda$ is a limit ordinal and $\lambda =
\lambda^{|\delta|}$ \then \, {\rm Ps}$'(\lambda,\delta)$. 

\noindent
4) If $\kappa = { \text{\rm cf\/}}(\delta),\kappa < \delta < \lambda,\lambda =
\lambda^\kappa$ and {\rm Pr}$_3(\lambda,\lambda,\lambda,\alpha)$ for
every $\alpha < \delta$, \then \, {\rm Ps}$'(\lambda,\delta)$. 

\noindent
5) {\rm Pr}$(\lambda,\kappa) \Rightarrow  { \text{\rm Ps\/}}(\lambda,\kappa)$,
{\rm Pr}$'(\lambda,\kappa) \Rightarrow { \text{\rm Ps\/}}'
(\lambda,\kappa)$ and similarly with $\delta$ instead of $\kappa$. 

\noindent
6) If $\lambda \ge 2^\kappa$ \then \,
$\lambda = \bold U_\kappa(\lambda) 
\Rightarrow { \text{\rm Pr\/}}(\lambda,\kappa)$.

\noindent
7) If {\rm Ps}$'(\lambda,\kappa)$ then {\rm Ps}$(\lambda,\kappa)$.
\end{claim}

\begin{remark}  Recall $\bold U_\kappa(\lambda) = \bold
U_{J^{\text{\rm bd}}_\kappa}(\lambda)$, and for an ideal $J$ on $\kappa,\bold
U_J(\lambda) = \text{ Min}\{|{\cP}|:{\cP} \subseteq
[\lambda]^\kappa$ is such that for every $f \in {}^\kappa \lambda$ for
some $A \in {\cP}$ we have $\{i < \kappa:f(i) \in A\} \ne \emptyset
\text{ mod } J\}$.
\end{remark}

\begin{PROOF}{\ref{2.5}}  
1) Let $\langle \lambda_i:
i < \text{ cf}(\lambda) \rangle$ be increasing
continuous with limit $\lambda$ such that $\delta^* < \lambda_0$ and
$2^{\lambda_i} < \lambda_{i+1}$, hence for limit
$\delta,\lambda_\delta$ is strong limit cardinal of cofinality cf$(\delta)$.
For $\delta \in S^{\text{cf}(\lambda)}
_{\text{cf}(\kappa)} = \{\delta < \text{ cf}(\lambda):\text{cf}(\delta) =
\text{ cf}(\kappa)\}$, let $\langle f_{\delta,\alpha}:\alpha < 
2^{\lambda_\delta} \rangle$ list the partial one-to-one
functions from $\lambda_\delta$
to $\lambda_\delta$ with domain of cardinality $\lambda_\delta$.  We choose
by induction on $\alpha < 2^{\lambda_\delta}$ a subset $A_{\delta,\alpha}$
of Dom$(f_{\delta,\alpha})$ of order type $\delta^*$ unbound in 
$\lambda_\delta$ such that $\beta < \alpha \Rightarrow \text{ sup}
(A_{\delta,\alpha} \cap A_{\delta,\beta}) < \lambda_\delta$; possible as
we have a tree with cf$(\delta)$ levels and $2^{\lambda_\delta}$ 
cf$(\delta)$-branches, each giving a possible $A_{\delta,\alpha}
\subseteq \text{ Dom}(f_{\delta,\alpha})$ and each
$A_{\delta,\beta} (\beta < \alpha)$ disqualifies $\le \lambda_\delta +
|\alpha|$ of them.

Now ${\cF} = \{f_{\delta,\alpha} \rest A_{\delta,\alpha}:
\delta \in S^{\text{cf}(\lambda)}
_{\text{cf}(\kappa)} \text{ and }\alpha < 2^{\lambda_\delta}\}$ is as
required because if $f$ is a partial function from $\lambda$
to $\lambda$ such that $|\text{Dom}(f)| = \lambda$ and $f$ is one to one then
$\{\delta < \text{ cf}(\lambda):(\exists^{\lambda_\delta} \, i <
\lambda_\delta)(i \in \text{ Dom}(f) \wedge 
f(i) < \lambda_\delta)\}$ contains a club of cf$(\lambda)$.

\noindent
2) This holds as $\diamondsuit_S$ for every stationary $S \subseteq \{\delta <
\lambda:\text{cf}(\delta) \ne \text{ cf}(\mu)\}$, see 
\cite{Sh:108}, and without any extra assumption by \cite{Sh:922}. 

\noindent
3) By the simple black box (see proof of \ref{1.4}, well it was phrased for
$\kappa$ but the same proof, and we have to rename $\lambda,\lambda^{<
\kappa}$ as $\lambda$; see \cite[Ch.IV]{Sh:e}, i.e. \cite{Sh:309}). 

\noindent
4) We combine the proof of the simple black box and the definition of Pr$_3$.
Let $\langle \gamma^*_i:i \le \kappa \rangle$ be increasing
continuous, $\gamma_0 = 0,\gamma_\kappa = \delta$.  By
Pr$_3(\lambda,\lambda,\lambda,\alpha)$ with $\alpha = \gamma^*_{i+1} -
\gamma^*_i$, for each $i < \kappa$, we can find ${\cF}_i$ such that
\mn
\begin{enumerate}
\item[$(*)_1$]   ${\cF}_i \subseteq \{g:g$ a partial function
from $\lambda$ to $\lambda$, otp$(\text{Dom}(g)) = \gamma^*_{i+1} -
\gamma^*_i\}$ 
\sn
\item[$(*)_2$]  $|{\cF}_i| \le \lambda$ 
\sn
\item[$(*)_3$]  for every $g^* \in {}^\lambda \lambda$ there
is $g \subseteq g^*$ from ${\cF}_i$
\end{enumerate}
\mn
By easy manipulation
\mn
\begin{enumerate}
\item[$(*)^+_3$]   if $g^* \in {}^\lambda \lambda$ and $\alpha <
\lambda$ \then \, there is $g \in {\cF}_i$ such that $g \subseteq g^*$
and Dom$(g) \subseteq [\alpha,\lambda)$.
\end{enumerate}
\mn
Clearly ${\cF}_i$ exists by the assumption
``Pr$_3(\lambda,\lambda,\lambda,\alpha)$ for 
$\alpha < \delta$" so let ${\cF}_i =
\{g_{i,\varepsilon}:\varepsilon < \lambda\}$.  Now for every $\eta \in
{}^\kappa\lambda$ let $f^0_\eta$ be the following partial function
from $({}^{\kappa >}\lambda) \times \lambda$ to $\lambda$:
\mn
\begin{enumerate}
\item[$(*)_4$]  if $i < \kappa,\varepsilon < \lambda$ and $\alpha \in
\text{ Dom}(g_{i,\eta(i)})$ \then \, $f^0_\eta((\eta \restriction
i,\alpha)) = g_{i,\eta(i)}(\alpha)$.
\end{enumerate}
\mn
Let $h$ be a one to one function from $({}^{\kappa >}\lambda) \times
\lambda$ onto $\lambda$ such that (if cf$(\lambda) \ge \delta$ then
also in $(*)_5(b)$ we can replace $\ne$ by $<$)
\mn
\begin{enumerate}
\item[$(*)_5$]   $(a) \quad \eta 
\in {}^{\kappa >}\lambda \wedge \alpha < \beta <
\lambda \Rightarrow h((\eta,\alpha)) < h(\eta,\beta)$ 
\sn
\item[${{}}$]  $(b) \quad \eta \triangleleft \nu \in {}^{\kappa >} \lambda
\wedge \alpha < \lambda \wedge \beta < \lambda \wedge \alpha \in
\text{ Dom}(g_{\ell g(\eta),\nu(\ell g(\eta))}) \Rightarrow
h((\eta,\alpha)) \ne$

\hskip25pt $h((\nu,\beta))$.
\end{enumerate}
\mn
Let $f_\eta$ be the following partial function from $\lambda$ to
$\lambda$ satisfying $f_\eta(\alpha) = f^0_\eta(h^{-1}(\alpha))$ so it
suffices to prove that $\cF = \{f_\eta:\eta \in {}^\kappa \lambda$ and
$f_\eta$ is one-to-one$\}$ exemplifies Ps$'(\lambda,\delta)$.

First, clearly each $f_\eta$ is a partial function from $\lambda$ to
$\lambda$.  Also for each $i < \kappa$ and $\varepsilon < \lambda$ the
function $g_{i,\varepsilon}$ has domain of order type $\gamma^*_{i+1} -
\gamma^*_i$, hence by $(*)_5(a)$ also $\eta \in {}^\kappa \lambda \wedge
i < \kappa \Rightarrow \text{ Dom}(f_\eta \rest \{h(\eta \rest
i,\varepsilon):\varepsilon < \lambda\}$ has order type $\gamma^*_{i+1} -
\gamma^*_i$.  By $(*)_5(b)$ also Dom$(f_\eta)$ has order type $\delta =
\sum\limits_{i}(\gamma^*_{i+1} \backslash \gamma^*_i)$.

Now if $f \in \cF$ then $f_\eta$ is one-to-one by the choice of $\cF$.
Second, let $f:\lambda \rightarrow \lambda$ be a one-to-one function
and we shall prove that for some $\eta,f_\eta \in \cF \wedge f_\eta
\subseteq f$.  We choose $\nu_i \in {}^i \lambda$ by induction on $i
\le \kappa$ such that $j<i \Rightarrow \nu_j = \nu_i \rest j$ and
$\nu_i \triangleleft \eta \in {}^\kappa \lambda \Rightarrow f_\eta
\rest \{h(\nu_j,\alpha):j<i,\alpha \in \text{ Dom}(g_{j,\nu_i(j)})\}
\subseteq f$.

For $i=0$ and $i$ limit this is obvious and for $i=j+1$ use
$(*)^+_3$.  So $\eta_\kappa \in {}^\kappa \lambda$ and
$f_{\eta_\kappa} \subseteq f$ hence is one-to-one hence
$f_{\eta_\kappa} \in \cF$ so we are done.

\noindent
5) Easy (recalling \ref{1.2A}(3)).

\noindent
6) Easy.

\noindent
7) Easy because if $f:\kappa \rightarrow \lambda$ is one-to-one then
   for some $u \subseteq \kappa$ of order type $\kappa,f \rest u$ is
   increasing (trivial if $\kappa$ is regular, easy if $\lambda$ is singular).
\end{PROOF}

\begin{claim}
\label{2.5A}  
1) Each of the following is a sufficient
condition to {\rm Pr}$_3(\chi,\lambda,\mu,\alpha)$, recalling
Definition \ref{2.4A}(2):
\mn
\begin{enumerate}
\item[$(a)$]   $\lambda^{|\alpha|} = \lambda = \chi \ge \mu > \alpha$
\sn
\item[$(b)$]   $\chi = \lambda \ge \mu > |\alpha|$ and $(\forall
\lambda_1 < \lambda)(\lambda^{|\alpha|}_1 < \lambda)$ 
\sn
\item[$(c)$]   $\chi = \lambda \ge \mu \ge \beth_\omega(|\alpha|)$.
\end{enumerate}
\mn
2) If $\chi_1 \le \chi_2,\lambda_1 \ge \lambda_2,\mu_1 \ge
\mu_2,\alpha_1 \ge \alpha_2$ \then \,
{\rm Pr}$_3(\chi_1,\lambda_1,\mu_1,\alpha_1)$ implies 
{\rm Pr}$_3(\chi_2,\lambda_2,\mu_2,\alpha_2)$.
\end{claim}

\begin{PROOF}{\ref{2.5A}}  
1) If clause $(a)$ holds, this is trivial, just use ${\cF}
= \{f:f$ a partial function from $\mu$ to $\lambda$ with $\alpha =
\text{ otp}(\text{Dom}(f))\}$.  If clause
$(b)$ holds, note that for every $f \in {}^\mu \lambda$, for some
$i_1,i_2 < \lambda$ we have $\alpha \le \text{ otp}(\{j < i_1:j < \mu$
and $f(j) < i_2\})$ and let ${\cF} = 
\{f:f$ a partial function from $\mu$
to $\lambda$ with bounded range and bounded domain if $\mu = \lambda$ 
such that $\alpha = \text{ otp}(\text{Dom}(f))\}$.
If clause $(c)$ holds, use \cite{Sh:460}. 

\noindent
2) Trivial.  
\end{PROOF}

\begin{claim}
\label{2.5B}  
1) In \ref{1.3}, \ref{1.4}(2) and in \ref{2.2} we
can weaken the assumption {\rm Pr}$(\lambda,\kappa)$ to
{\rm Ps}$(\lambda,\kappa)$. 

\noindent
2) In \ref{2.4} we can weaken the assumption {\rm Pr}$'(\lambda,\delta^*)$
to {\rm Ps}$'(\lambda,\delta^*)$.
\end{claim}

\begin{PROOF}{\ref{2.5B}}  
The same proofs.  

We can get another answer on the existence of universals.
\end{PROOF}

\begin{claim}
\label{2.5C}  
If $\lambda$ is strong limit, {\rm cf}$(\lambda) > 
{ \text{\rm cf\/}}(\kappa)$ and 
$\lambda > \theta(\ge \kappa)$, \then \, in ${\gH}_{\lambda,
\theta,\kappa}$ there is no {\rm we}-universal member even for the
class of bipartite members. 
\end{claim}

\begin{PROOF}{\ref{2.5C}}  Let $\delta := \theta^+ \times
\kappa$ (recalling $\lambda$ is a limit cardinal) by 
\ref{2.5}(1) we have Pr$'(\lambda,\delta)$
hence by \ref{2.4} we are done.  
\end{PROOF}

\noindent
Note that Ps may fail.
\begin{claim}
Assume $\delta < \lambda$, {\rm cf}$(\lambda) \le \text{\rm
cf}(\delta)$ and $\alpha < \lambda \Rightarrow |\alpha|^{\text{\rm
cf}(\delta)} < \lambda$.

\Then \, {\rm Ps}$(\lambda,\delta)$ fail (hence also {\rm
Pr}$(\lambda,\delta)$, {\rm Pr}$'(\lambda,\delta)$, {\rm
Ps}$'(\lambda,\delta)$. 
\end{claim}

\begin{proof}  Toward contradiction let $\cF$ witness
Pr$(\lambda,\delta)$ so $\cF$ is of cardinality $\lambda$.  Let
$\langle f_\varepsilon:\varepsilon < \lambda\rangle$ list $\cF$;
choose an increasing sequence $\langle \lambda_i:i < \text{
cf}(\lambda)\rangle$ such that $\lambda_i =
(\lambda_i)^{\text{cf}(\delta)}$ and $\lambda =
\varepsilon\{\lambda_i:i < \text{cf}(\lambda)\}$.  We choose $\cU_i$
by induction on $i$ such that:
\mn
\begin{enumerate}
\item[$(*)^1_i$]  $(a) \quad \cU_i \subseteq \lambda$
\sn
\item[${{}}$]  $(b) \quad \lambda_i \subseteq \cU_i$
\sn
\item[${{}}$]  $(c) \quad |\cU_i| = \lambda_i$
\sn
\item[${{}}$]  $(d) \quad$ if $f \in \cF$ and Dom$(f) \cap \cU_i$ is
unbounded in Dom$(f)$, \then \, Dom$(f) \cup \text{ Rang}(f) \subseteq 
\lambda$.
\end{enumerate}
\mn
For clause (d) note that if Dom$(f) \cap \cU_i$ is unbounded in
Dom$(f)$ then there is $u \subseteq \text{ Dom}(f) \cap \cU_i$
unbounded in Dom$(f)$ of order-type cf$(\gamma)$ and such $u$
determines $f$ in $\cF$ uniquely

Now choose $f^*:\lambda \rightarrow \lambda$ such that $f^*$ maps
$\lambda[\bigcup\limits_{j < i(*)},\lambda_{i(*)})(\cU_i \backslash
\cup\{\cU_j:j < i\})$ into $[\lambda_{i(1)},\lambda_{i(j)+1})$ when
$i(2) \le i(j) < \text{\rm cf}(\lambda)$ and be increasing $f^*$
contradict the choice of $\cF$.
\end{proof}

\section {Complete characterization under GCH} 

We first resolve the case $\lambda$ is strong limit and get a complete
answer in \ref{2.9} by dividing to cases (in \ref{2.6},
\ref{2.7}, \ref{2.8} and \ref{2.5C}), in \ref{2.8} we deal
also with other cardinals.  This includes cases in which there are
universals (\ref{2.6}, \ref{2.7}) and the existence of
we-universal and of ste-universal are equivalent.  In fact in
\ref{2.8} we deal also (in part (2)) with another case:
$\kappa = \theta \le \text{ cf}(\lambda),\lambda =
\sum\limits_{\alpha < \lambda,\beta < \theta} |\alpha|^{|\beta|}$ (and
then there is no we-universal).

Next we prepare the ground for resolving the successor case under GCH
(or weaken conditions using also \ref{2.8}(2)).  If $\lambda = \mu^+
= \theta,\mu = \mu^\kappa$ there is a we-universal in
$\gH_{\lambda,\theta,\kappa}$ (\ref{2.10}), if
$\lambda = \mu^+ = 2^\mu,\kappa < \mu$ (in \ref{2.11}, \ref{2.12})
we give a sufficient condition for existence.  In \ref{2.14} we
sum up.  We end with stating the conclusion for the classes of bipartite
graphs (\ref{2.15}, \ref{2.16}).

\begin{claim}
\label{2.6}  
Assume $\lambda$ is strong limit, $\text{\rm cf}(\lambda) \le
\text{\rm cf}(\kappa),\kappa \le \theta < \lambda$ and
$\kappa < \theta \vee \text{\rm cf}(\lambda) < \text{\rm cf}(\kappa)$;
hence {\rm cf}$(\lambda) < \theta < \lambda$ so $\lambda$ is singular. 

\noindent
\Then \, in ${\gH}_{\lambda,\theta,\kappa}$ there is a
{\rm ste}-universal member.
\end{claim}

\begin{PROOF}{\ref{2.6}}
  Denote $\sigma = \text{ cf}(\lambda)$ and let $\langle
\lambda_i:i < \sigma \rangle$ be increasing continuous with limit $\lambda$
such that $\lambda_0 > \theta$ and $(\lambda_{i+1})^{\lambda_i} = 
\lambda_{i+1}$.  For any graph $G \in {\gH}_{\lambda,\theta,\kappa}$ we
can find $\langle V^G_i:i < \sigma \rangle$ such that:  $V^G_i \subseteq
V^G,\langle V^G_i:i < \sigma\rangle$ is 
increasing continuous with $i$ with union $V^G$ such that 
$|V^G_i| = \lambda_i$ and
\mn
\begin{enumerate}
\item[$(*)_1$]   if $x \in V^G \backslash V^G_{i+1}$ then
$|\{y \in V^G_{i+1}:y$ is $G$-connected to $x\}| < \kappa$. 
\end{enumerate}
\mn
As $\sigma = \text{ cf}(\lambda) \le \text{ cf}(\kappa)$ it follows
\mn
\begin{enumerate}
\item[$(*)_2$]   if $i < \sigma$ is a limit ordinal and $x \in V^G 
\backslash V^G_i$ \then \, (cf$(i) < \sigma \le \text{ cf}(\kappa)$ hence)
$|\{y \in V^G_i:y$ is $G$-connected to $x\}| < \kappa$.
\end{enumerate}
\mn
For $i \le \sigma$ let 
$T_i = \prod\limits_{j < i} 2^{\lambda_{j+1}},T^s = \cup\{ T_{i+1}:
i < \sigma\},T = \bigcup_{i < \sigma} T_i$.

Let $\bar A = \langle A_\eta:\eta \in T^s \rangle$ be a sequence of pairwise
disjoint sets such that $\eta \in T_{i+1} \Rightarrow |A_\eta| =
\lambda_i$.  
For $\eta \in T \cup T_\sigma$
let $B_\eta = \cup \{A_{\eta \restriction j}:j < \ell g
(\eta),j$ a successor ordinal$\}$.  Now we choose by induction on 
$i \le \sigma$, for each $\eta \in T_i$ a graph $G_\eta$ such that:
\mn
\begin{enumerate}
\item[$\boxplus$]   $(a) \quad V^{G_\eta} = B_\eta$ 
(so for $\eta = <>$ this is the graph with the empty set of

\hskip25pt  nodes) and so $|V^{G_\eta}| = \Sigma\{\lambda_j:j < \ell g(\eta)$

\hskip25pt  $\text{ successor}\}$
\sn
\item[${{}}$]  $(b) \quad$ if 
$\nu \triangleleft \eta$ then $G_\nu$ is an induced
subgraph of $G_\eta$, moreover 

\hskip25pt $(\forall x \in V^{G_\eta} \backslash V^{G_v})(x \text{ is } 
G_\eta \text{-connected to } < \kappa$ nodes in $V^{G_\nu})$
\sn
\item[${{}}$]  $(c) \quad$ if $i < \sigma,\eta \in T_i,G$ a graph such that
$|V^G| = \lambda_{i+1}$ and $G \in \gH_{\lambda_{i+1},\theta,\kappa}$

\hskip25pt  and $G_\eta$ is an induced subgraph of $G$ and
$(\forall x \in V^G \backslash V^{G_\eta})(x$ is 

\hskip25pt $G$-connected to $< \kappa$
members of $V^{G_\eta})$, \then \, for some $\alpha < 2^{\lambda_{i+1}}$

\hskip25pt there is  an isomorphism from $G$ 
onto $G_{\eta \char 94 \langle \alpha
\rangle}$ which is the 

\hskip25pt identity on $B_\eta$
\sn
\item[${{}}$]  $(d) \quad G_\eta \in {\gH}_{|G_\eta|,\theta,\kappa}$. 
\end{enumerate}
\mn
[Can we carry the induction?  For $i=0$ this is trivial.  For $i=j+1$
this is easy, the demand in clause (c) poses no threat to the others.
For $i$ limit for $\eta \in T_i$, the graph $G_\eta$ is well defined
satisfying clauses (a), (b) (and (c) is irrelevant), but why $G_\eta
\in {\frak H}_{|G_\eta|,\theta,\kappa}$?  Toward contradiction assume
$A_0,A_1 \subseteq B_\eta,A_0 \times A_1 \subseteq R^{G_\eta}$ and
$\{|A_0|,|A_1|\} = \{\kappa,\theta\}$.  If $\ell < 2$ and
cf$(|A_\ell|) \ne \text{ cf}(i)$ then for some $j <i,|A_\ell \cap
B_{\eta \restriction j}| = |A_\ell|$ so \wilog \, $A_\ell \subseteq
B_{\eta \restriction j}$, but then by clause (b) no 
$x \in B_\eta \backslash B_{\eta \restriction j}$ 
is $G_\eta$-connected to $\ge \kappa$
members of $B_{\eta \restriction j}$ and $|A_\ell| \ge 
\text{ Min}\{\kappa,\theta\} = \kappa$, hence $A_{1-\ell} \subseteq B_{\eta
\restriction j}$, so we get contradiction to the induction hypothesis.
So the remaining case is cf$(|A_0|) = \text{ cf}(i) = \text{
cf}(|A_1|)$ hence cf$(\theta) = \text{ cf}(\kappa) = \text{ cf}(i)$;
so as we are assuming cf$(\kappa) \ge \text{ cf}(\lambda) 
\ge \text{ cf}(i)$, we
necessarily get cf$(i) = \text{ cf}(\lambda) = \text{ cf}(\kappa) = \text{
cf}(\theta)$.  By the last assumption of the claim (i.e. $\kappa < \theta \vee
\text{ cf}(\kappa) > \text{ cf}(\lambda))$ we get that $\kappa <
\theta$ and \wilog \, $|A_0| = \kappa,|A_1| = \theta$, so for some $j
< \text{ cf}(\lambda)$ we have $j < i$ and 
$|A_1 \cap B_{\eta \restriction j}| \ge
\kappa$, so as above $A_0 \subseteq B_{\eta \restriction j}$, hence
again as above $A_1 \subseteq B_{\eta \restriction j}$ and we are done.]
\medskip

\noindent
We let $G^* = \cup\{G_\eta:\eta \in T\}$.  Now we
shall check that $G^*$ is as required.  First assume toward
contradiction that $A,B \subseteq V^{G^*}$
and $A \times B \subseteq R^{G^*}$ and $\{|A|,|B|\} = \{\kappa,\theta\}$.
A set $C \subseteq V^{G^*}$ will be called $\bar A$-flat if it is
included in some $B_\eta,\eta \in T \cup T_\sigma$.  
Easily above if $B$ is not $\bar A$-flat then $A$ is
$\bar A$-flat.  So \wilog \, for some $\eta \in T_\sigma$ we have
$A \subseteq B_\eta$ but then
$x \in G^* \backslash B_\eta \Rightarrow$ (for some $i(x) < \sigma$
and $\nu \in T_{i(x)}$ we have $x \in B_\nu \backslash B_{\eta \cap \nu})
\Rightarrow |\{y \in B_\eta:x$ is
$G^*$-connected to $y\}| < \kappa$, so as $\kappa \le \theta$ we get
$B \subseteq B_\eta$, hence $A,B \subseteq V^{G_\eta}$ 
and we get contradiction to clause (d) of $\boxplus$.   

\noindent
So $K_{\kappa,\theta}$ does
not weakly embed into $G^*$; also $|V^{G^*}| = \lambda$ so $G^* \in
{\gH}_{\lambda,\theta,\kappa}$.  Lastly, the 
ste-universality follows from
the choice of $\langle V^G_i:i < \sigma \rangle$ for any $G \in
{\gH}_{\lambda,\theta,\kappa}$.  That is, we can choose by induction on
$i,\eta_i \in T_i$ and an isomorphism $f_i$ from $G \restriction V^G_i$
onto $G_{\eta_i}$ if $i > 0$, with $f_i$ increasing continuous (and $\eta_i$
increasing continuous) using for successor $i$ clause (c) of $\boxplus$.
\end{PROOF}

In the previous claim we dealt with the case of $\kappa,\theta <
\lambda$.  In the following claim we cover the case of $\theta = \lambda$:
\begin{claim}
\label{2.7}  
Assume $\lambda$ is strong limit, {\rm cf}$(\lambda) 
\le { \text{\rm cf\/}}(\kappa) \le \kappa < \theta = \lambda$; hence 
$\lambda$ is singular.  

\noindent
\Then \, in
${\gH}_{\lambda,\theta,\kappa}$ there is a {\rm ste}-universal member.
\end{claim}

\begin{PROOF}{\ref{2.7}}  
Similar to the previous proof and
\cite[Th.3.2, p.268]{Sh:26}.  
Let $\sigma = \text{ cf}(\lambda)$ and 
$\bar \lambda = \langle \lambda_i:i < \sigma \rangle,\langle T_i:i \le
\sigma\rangle,T^s,T$ be as in the proof of \ref{2.6}.
For any graph $G \in {\gH}_{\lambda,\theta,\kappa}$ let $h^G:{}^\kappa(V^G)
\rightarrow \sigma$ be defined by: if $|\{x_\varepsilon:\varepsilon <
\kappa\}| = \kappa$ then $h^G(\bar x)$ is the first $i < \sigma$
such that $\lambda_i \ge |\{y \in V^G:y$ is $G$-connected to every $x_i,
i < \kappa = \ell g(\bar x)\}|$, otherwise $h^G(\bar x)$ is not
defined.  Now we choose $\langle V^G_i:i < \sigma \rangle$
as an increasing continuous sequence of subsets of $V^G$ with union $V^G$
such that if $i < \sigma$ then
$|V^G_i| \le \lambda_i$ and $\bar x \in {}^\kappa(V^G_{i+1}) \wedge
|\text{Rang}(\bar x)| = \kappa \wedge 
h^G(\bar x) \le i +1 \Rightarrow (\forall y \in V^G)
(``y$ is $G$-connected to $x_i$ for every $i < \kappa" \Rightarrow y \in
V^G_{i+1})$.

Then when (as in the proof of \ref{2.6}) we construct $\langle G_\eta:
\eta \in T^s \rangle$ we also construct
$\langle h_\eta:\eta \in T^s \rangle,h_\eta:{}^\kappa(B_\eta) \rightarrow
\sigma$ with the natural demands.  In the end we have to check
that ``$K_{\kappa,\theta}$ is not strongly embeddable into 
$G^*$"; if cf$(\kappa) = \sigma$ we need to look at
slightly more (as in the end of the proof of \ref{2.6}).
\end{PROOF}

\begin{remark}  More generally see \cite{Sh:829}.
\end{remark}

\begin{claim}
\label{2.8}  
1) Assume $\lambda$ is strong limit,
$\lambda \ge \theta = \kappa$, {\rm cf}$(\kappa) = 
{ \text{\rm cf\/}}(\lambda)$, \then \, in
${\gH}_{\lambda,\theta,\kappa}$ there is no {\rm we}-universal graph,
even for the members from ${\gH}^{\text{\rm sbp}}_{\lambda,\theta,\kappa}$. 

\noindent
2) Assume
\mn
\begin{enumerate}
\item[$(a)$]  $\kappa = \theta \le \lambda$
\sn
\item[$(b)$]   $(\forall \alpha < \lambda)(\forall \beta < \theta)
[|\alpha|^{|\beta|} \le \lambda]$, (recall $\kappa = \theta$)
\sn
\item[$(c)$]   {\rm cf}$(\lambda) \ge { \text{\rm cf\/}}(\kappa)$.
\end{enumerate}
\mn
\Then \, in ${\gH}_{\lambda,\theta,\kappa}$ there is no
we-universal graph, even for the 
${\gH}^{\text{\rm sbp}}_{\lambda,\theta,\kappa}$.
\end{claim}

\begin{PROOF}{\ref{2.8}}
  1) By part (2).

\noindent
2) Let $\sigma = \text{ cf}(\kappa)$.
Let $\langle \gamma_i:i < \sigma \rangle$ be (strictly)
increasing with limit $\theta = \kappa$.  

Without loss of generality
\mn
\begin{enumerate}
\item[$\boxtimes$]   $|\gamma_{i+1} - \gamma_i|$ is (finite or) a
cardinal $(< \kappa)$ with cofinality $\ne \text{ cf}(\lambda)$.
\end{enumerate}
\mn
[Why?  If $\kappa$ is a limit cardinal, trivial as if a cardinal 
$<\kappa$ fails its successor is O.K.; if $\kappa$ is a
successor cardinal, $\gamma_i = i$ is
O.K. and also $\gamma_i = \omega i$ or $\gamma_i = \omega_1 i$ is O.K.]

Given a graph $G^*$ in ${\gH}_{\lambda,\theta,\kappa}$ \wilog \,
$V^{G^*} = \lambda$.

For $i < \sigma$ let $T_i = \{f:f$ is a partial one-to-one mapping 
from $\gamma_i$ into $\lambda = V^{G^*}$ with bounded range such that 
$j<i \wedge \ell \in\{0,1\} \Rightarrow |\gamma_{j+1} - \gamma_j|/2=1,
\gamma_j \le \varepsilon \le
\gamma_{j+1}$ and $\varepsilon - \ell$ mod $2\}|$ and $2 \alpha,
2 \beta +1 \in \text{ Dom}(f)
\Rightarrow f(2 \alpha) R^{G^*} f(2 \beta +1)\}$.  Let 
$T = \bigcup\limits_{i < \sigma} T_i$, so $T$ is a 
tree with $\le \lambda$ nodes and $\sigma$ levels; for $\eta \in T$
let $\bold i(\eta) < \sigma$ be the unique $i < \sigma$ such that
$\eta \in T_i$.
Let $T^+ = \{\eta:\text{for some } \zeta,\ell g(\eta) = \zeta +1,\eta
\restriction \zeta \in T$ is $\triangleleft$-maximal in $T$ and $\eta(\zeta)
= 0\}$.  

Note
\mn
\begin{enumerate}
\item[$(*)_1$]   if $i < \sigma$ is a limit ordinal, $\langle f_j:j
< i \rangle$ is $\subseteq$-increasing, $f_j \in T_j$ 
then $\bigcup\limits_{j < i} f_j \in
T_i$ [in other words if $f$ is a function from $\gamma_i$ to $\lambda$
such that $j < i \Rightarrow f \restriction \gamma_j \in T_j$ then
$f \in T_i$]. 
\end{enumerate}
\mn
[Why?  The least obvious demand is sup Rang$(f) < \lambda$ which
holds as cf$(i) \le i < \sigma = \text{ cf}(\kappa) \le \text{ cf}(\lambda)$.]
\mn
\begin{enumerate}
\item[$(*)_2$]   there is no $\bar f = \langle f_i:i < \sigma
\rangle$ increasing such that $i < \sigma \Rightarrow f_i \in T_i$.
\end{enumerate}
\mn
[Why?  As then $\bigcup\limits_{i < \sigma} f_i$ weakly embed a complete
$(\kappa,\theta)$-bipartite graph into $G^*$.]

We define a bipartite graph $G$

\begin{equation*}
\begin{array}{clcr}
U^G = &\{\eta:\eta \in T,\bold i(\eta) \text{ is even}\} \\
  &\cup \{(\eta,\varepsilon):\eta \in T^+,\varepsilon < 
\lambda \text{ is even}\}
\end{array}
\end{equation*}

\begin{equation*}
\begin{array}{clcr}
V^G = &\{\eta:\eta \in T,\bold i(\eta) \text{ is odd}\} \\
  &\cup \{(\eta,\varepsilon):\eta \in T^+ \text{ and } \varepsilon <
\lambda \text{ is odd}\}
\end{array}
\end{equation*}

\begin{equation*}
\begin{array}{clcr}
R^G = R^G_1 \cup R^G_2 \text{ where } R^G_1 = 
\bigl\{ \{\nu,(\eta,\varepsilon)\}:&\eta \in T^+,\bold i(\eta) 
\text{ is a succcessor ordinal and } \nu \in U^G \cap T \\
  &\text{ (and } \gamma_{\bold i(\eta)-1} \le \varepsilon <
\gamma_{\bold i(\eta)}) \text{ and} \\
  &\nu \triangleleft \eta \text{ and } 
\varepsilon \ne \ell g(\nu) \text{ mod } 2 \bigr\} \\
  &F^G_2 = \bigl\{\{(\eta,\varepsilon_1),(\eta,\varepsilon_2)\}:
\eta \in T^+,\bold i(\eta) \text{ is a successor ordinal,} \\
  &\{\varepsilon_1,\varepsilon_2\} \subseteq [\bold i(\eta)-1,\bold
i(\eta) \text{ and } \varepsilon_1 \ne \varepsilon_2 \text{ mod }2 \bigr\}.
\end{array}
\end{equation*}

Now
\mn
\begin{enumerate}
\item[$(a)$]   $|T| \le \lambda$ by clause (b) of the assumption,
$[\gamma_{\bold i(\eta)-1},\gamma_{\bold i(\eta)}]$
which is a weak form of ``$\lambda$ is strong limit"
\sn
\item[$(b)$]    $|T^+| \le \lambda$.
\end{enumerate}
\mn
If $(\cU_1 \cap (T^+ \times \lambda)) \ne \emptyset \ne (V_1 \cap (T^+
\times \lambda))$ then we can choose $(\eta_1,\varepsilon_1) \in
\cU_,(\eta_2,\varepsilon_2) \in V_1$.
\mn
[Why?   As $\eta \in T^+ \Rightarrow \text{ sup Rang}(\eta) < 
\lambda$, see $(*)_2$ recalling $i = \ell g(\eta) < \sigma$,
Rang$(\eta) \subseteq \gamma_i < \lambda$.]
\mn
\begin{enumerate}
\item[$(c)$]    $|U^G| \ge |T_1| = \lambda$ and $|V^G| \ge |T_2| = \lambda$
\end{enumerate}
\mn
hence
\mn
\begin{enumerate}
\item[$(d)$]   $|U^G| = \lambda = |V^G|$ so
\sn
\item[$(e)$]   $G \in {\gH}^{\text{\rm sbp}}_{\lambda,\theta,\kappa}$.
\end{enumerate}
\mn
[Why?  Being bipartite is obvious; so toward contradiction assume $(U_1
\subseteq U^G \wedge V_1 \subseteq V^G) \vee (U_1 \subseteq V^G
\wedge V_1 \subseteq U^G)$ have cardinality $\kappa$ (recall that
$\kappa = \theta$) and $U_1 \times V_1
\subseteq R^G$.  If $(\eta_*,\varepsilon_*) \in V_1$ 
and $i < \text{\rm cf}(\kappa) =
\sigma$ be such that $\varepsilon_* < \gamma_i$, then $U_1 \subseteq
\{\nu:\nu \triangleleft \eta\} \cup \{(\nu,\zeta) \in U^G \cup V^G:\nu
= \eta_*$ and $\zeta < \gamma_i\}$ which clearly has cardinality $<
\kappa$, contradiction.  Hence $V_1 \cap (T^+ \times \lambda) =
\emptyset$ and by symmetry $\cU_1 \cap (T^+ \cap \lambda) =
\emptyset$, hence $V_1 \subseteq T,U_1 \subseteq T$ hence $U_1 \times
V_1$ (set of unordered pairs) is disjoint to $R$, contradiction.] 
\mn
\begin{enumerate}
\item[$(f)$]   $G$ is not weakly embeddable into $G^*$.
\end{enumerate}
\mn
[Why?  If $f$ is such an embedding, we try to choose by induction on
$i < \sigma$, a member $\eta_i$ of $T_i$, increasing continuous with
$i$ such that $(\forall \varepsilon \in \text{ Dom}(\eta_i))(\forall j
< i)[\gamma_j \le \varepsilon < \gamma_{j+1} \Rightarrow 
\eta_i(\varepsilon) = f((\eta_i \restriction \gamma_j,\varepsilon))]$. 
If we succeed we get a contradiction to $G^* \in 
{\gH}_{\lambda,\theta,\kappa}$ by $(*)_2$, so we cannot carry the
induction for every $i < \sigma$.  For $i=0$ and 
$i$ limit there are no problems 
(see $(*)_1$), so for some $i=j+1 < \sigma,f_j$ is well defined but we
cannot choose $f_i$.  But if $j < \sigma$ consider $f_i = f_j \cup
\{(\varepsilon,f((\eta_j,\varepsilon)):\varepsilon \in
[\gamma_j,\gamma_i)\}$. This gives a contradiction except possibly
when $\lambda = \text{ sup Rang}(\eta_i)$, but then necessarily
by $\boxtimes$, $|\gamma_{i+1} - \gamma_i|$ has cofinality $\ne \text{
cf}(\lambda)$, so for $\ell < 2$, for some $\iota_\ell < \sigma$ the
set $\{\varepsilon:\gamma_j \le \varepsilon < \gamma_{j+1}$ and
$f(\eta_j,\varepsilon) < \lambda_{\iota_\ell}$ and $\varepsilon 
= \ell \text{ mod } 2\}$ has cardinality $|\gamma_i - \gamma_j|/2$
which has cofinality $\ne \text{ cf}(\lambda)$, and then $f_i := f_j \cup
\{(\varepsilon,f(\eta_j,\varepsilon)):\gamma_j \le \varepsilon <
\gamma_i$ and $f(\eta_j,\varepsilon) < 
\lambda_{\text{\rm max}\{\iota_0,\iota_1\}}\}$ is O.K., contradiction.]
 \end{PROOF}

\begin{theorem}
\label{2.9}  
Assume $\lambda \ge \theta \ge \kappa \ge
\aleph_0$ and $\lambda$ is a strong limit cardinal.  
There is a we-universal in ${\gH}_{\lambda,\theta,\kappa}$ \Iff \, 
{\rm cf}$(\lambda) \le \text{\rm cf}(\kappa)$ and 
$(\kappa < \theta \vee \text{\rm cf}(\lambda) <
\text{\rm cf}(\theta))$ \Iff \,  there is a ste-universal in
${\gH}_{\lambda,\theta,\kappa}$.
\end{theorem}

\begin{remark}
\label{2.9d}
Similarly for the universal for $\{G^{[\text{\rm
gr}]}:g \in \gH^{\text{\rm sbp}}_{\lambda,\theta,\kappa}\}$.
\end{remark}

\begin{PROOF}{\ref{2.9}}  
We use freely \ref{0.5}(3) and below in each case the
middle condition in \ref{2.9} clearly holds or clearly fails and the other
conditions hold or fail by the claim quoted in the case.

If $\theta < \lambda$ and $\text{cf}(\lambda) > \text{ cf}(\kappa)$ by 
\ref{2.5C} the family ${\gH}_{\lambda,\theta,\kappa}$ has no
we-universal.

If $\theta = \lambda$ and $\text{cf}(\lambda) > \text{ cf}(\kappa)$
then Pr$'(\lambda,\kappa)$ holds by \ref{2.5}(1) (and recall that
Pr$(\lambda,\kappa)$ is equivalent here to Pr$'(\lambda,\kappa)$,
since $\kappa$ is a cardinal, see \ref{1.2A}(1)(iii)) hence by
\ref{2.2} the family $\gH_{\lambda,\theta,\kappa}$ has no
we-universal member.

If cf$(\lambda) < \text{ cf}(\kappa)$ and
$\theta < \lambda$ by \ref{2.6} the family
$\gH_{\lambda,\theta,\kappa}$ has a ste-universal member; 
the second statement in Theorem \ref{2.9} holds as:
if $\kappa < \theta$ easy, if $\kappa \ge \theta$ then $\kappa =
\theta$ hence cf$(\lambda) < \text{ cf}(\kappa) = \text{ cf}(\theta)$.

If cf$(\lambda) < \text{ cf}(\kappa)$ (hence $\kappa \ne \lambda$ so
$\kappa < \lambda$) and $\theta = \lambda$ (so $\kappa < \theta$) 
by \ref{2.7} the family
$\gH_{\lambda,\theta,\kappa}$ has a ste-universal member.

So the remaining case is cf$(\lambda) = \text{ cf}(\kappa)$.  
If $\kappa < \theta < \lambda$ by \ref{2.6} in
$\gH_{\lambda,\theta,\kappa}$ there is an ste-universal; 
if $\kappa = \theta \le \lambda$ by
\ref{2.8}(1) in $\gH_{\lambda,\theta,\kappa}$ there is no we-universal
member.  If $\kappa < \theta = \lambda$ by \ref{2.7} in
$\gH_{\lambda,\theta,\kappa}$ there is a ste-universal member.  We
have checked all possibilities hence we are done.
\end{PROOF}

We turn to $\lambda$ successor cardinal.  In the following case,
possibly the existence of we-universal and ste-universal are not
equivalent, see \ref{2.12}(2) + \ref{2.12}(4) and Theorem \ref{2.14}.
\begin{claim}
\label{2.10}  
Assume $(\lambda \ge \theta \ge \kappa \ge \aleph_0$ and)
\mn
\begin{enumerate}
\item[$(a)$]  $\lambda = \mu^+$
\sn
\item[$(b)$]  $\kappa < \mu$ and $\theta = \lambda$
\sn
\item[$(c)$]  $\mu = \mu^\kappa$.
\end{enumerate}
\mn
\Then \, in ${\gH}_{\lambda,\theta,\kappa}$ there is a we-universal
member. 
\end{claim}

\begin{PROOF}{\ref{2.10}}
If $G \in {\gH}_{\lambda,\theta,\kappa}$ (and
\wilog \, $V^G = \lambda$) and $\alpha < \lambda$, then
\mn
\begin{enumerate}
\item[$(*)_\alpha$]  $\{\beta < \lambda:\beta$ is $G$-connected 
to $\ge \kappa$ elements $\gamma < \alpha\}$ is bounded in 
$\lambda$ say by $\beta_\alpha < \lambda$.
\end{enumerate}
\mn
Hence there is a club $C = C_G$ of $\lambda$ such that:
\mn
\begin{enumerate}
\item[$(i)$]   cf$(\alpha) \ne \text{ cf}(\kappa),
\alpha \in C,\beta \in [\alpha,
\lambda) \Rightarrow \kappa > \text{ otp}\{\gamma < \alpha:\gamma$ is
$G$-connected to $\beta\}$
\sn
\item[$(ii)$]  cf$(\alpha) = \text{ cf}(\kappa),
\alpha \in C,\beta \in [\alpha,
\lambda) \Rightarrow \kappa \ge \text{ otp}\{\gamma < \alpha:\gamma$ is
$G$-connected to $\beta\}$
\sn
\item[$(iii)$]  if $\alpha \in C$ then $\mu/\alpha$ and $\alpha >
\sup(C \cap \alpha) \Rightarrow \text{ cf}(\alpha) \ne \text{
cf}(\kappa)$.
\end{enumerate}
\mn
We shall define $G^*$ with $V^{G^*} = \lambda$ below.  For each $\delta <
\lambda$ divisible by $\mu$ let $\langle a^\delta_i:
i < \mu \rangle$ list ${\cP}_\delta = 
\{a:a \subseteq \delta$, and $|a| < \kappa$ or
otp$(a) = \kappa$ and $\delta = \sup(a)\}$, each appearing $\mu$ times,
possible as $|\delta| = \mu = \mu^\kappa$, and let

\begin{equation*}
\begin{array}{clcr}
R^{G^*}_\delta = \bigl\{ &\{\beta,\delta +i\}:\delta < \lambda \text{
is divisible by } \mu,\beta < \delta,i < \mu,\beta
\in a^\delta_i \bigr\} \\
  &\cup \bigl\{ \{\delta +i,\delta +j\}:i \ne j <
\mu \text{ and } \delta < \lambda \text{ is divisible by } \mu \bigr\}.
\end{array}
\end{equation*}
\mn
Now clearly we have $\alpha + \mu \le \beta < \lambda \Rightarrow
\kappa > |\{\gamma < \alpha:\gamma \text{ is } G^* \text{-connected to }
\beta\}|$ hence $K_{\kappa,\lambda}$ (which is $K_{\kappa,\theta}$ by
the assumptions) cannot be weakly embedded into $G^*$.  On
the other hand if $G \in {\gH}_{\lambda,\theta,\kappa}$ \wilog \,
$V^G = \lambda$ and let $C_G$ be as above, 
and let $\langle \alpha_\zeta:\zeta < \lambda \rangle$
list in increasing order $C_G \cup \{0\}$, and we can choose by induction on
$\zeta$, a weak embedding $f_\zeta$ of $G \restriction \alpha_\zeta$ into
$G^* \restriction (\mu \times \zeta)$.  So $G^*$ is as required. 
\end{PROOF}

\begin{claim}
\label{2.11}  
Assume $(\lambda \ge \theta \ge \kappa \ge \aleph_0$ and)
\mn
\begin{enumerate}
\item[$(a)$]   $\lambda = 2^\mu = \mu^+,\mu$ is a singular cardinal
\sn
\item[$(b)$]  $\kappa < \mu$ and \footnote{in fact, $\kappa = \theta$ is
O.K., but already covered by \ref{2.8}(2)} $\kappa < \theta \le \lambda$
\sn
\item[$(c)$]   for every ${\cP} \subseteq [\mu]^\mu$ of
cardinality $\lambda$ for some \footnote{note that if
$\beth_\omega(\kappa) \le \mu$ this clause always holds; and if
$2^\kappa \le \mu$ it is hard to fail it, not clear if its negation is
consistent} $B \in [\mu]^\kappa$, for $\lambda$ sets $A
\in {\cP}$ we have $B \subseteq A$
\sn
\item[$(d)$]   {\rm cf}$(\kappa) = { \text{\rm cf\/}}(\mu)$.
\end{enumerate}
\mn
\Then \, in ${\gH}_{\lambda,\theta,\kappa}$ there is no 
{\rm we}-universal member (even for the family of bipartite graphs).
\end{claim}

\begin{PROOF}{\ref{2.11}}  
Let $G^* \in {\gH}_{\lambda,\theta,\kappa}$, \wilog \,
$V^{G^*} = \lambda$ and we shall construct a $G \in 
{\gH}^{\text{sbp}}_{\langle \mu,\lambda \rangle,
\theta,\kappa}$ not weakly embeddable into it.
Now we choose $U^G = \mu,V^G = \lambda \backslash \mu$.

Notice that $\lambda^\mu = \lambda$ (by (a)), so let 
$\langle (f_\alpha,B_\alpha):\mu \le \alpha < \lambda \rangle$ list
the pairs $(f,B)$ such that $f:\mu \rightarrow \lambda$ is one to one,
$B \in [\mu]^\mu$ and $f \restriction B$ is increasing such that each pair
appears $\lambda$ times.  Let $\beta_B = \sup\{\beta +1:\beta <
\lambda$ is $G^*$-connected to $\mu$ members of $B\}$ for $B \in
[\lambda]^\mu$ (and we shall use it for $B \in [\lambda \backslash
\mu]^\mu$, i.e. for subsets of $V^G$).  Now $\beta_B$ 
is $< \lambda$ by clause (c) of the assumption.  We shall now choose
inductively $C_\alpha$ for $\alpha \in [\mu,\lambda)$ such that
\mn
\begin{enumerate}
\item[$\circledast$]   $(i) \quad C_\alpha \subseteq 
B_\alpha$ is unbounded of order type $\kappa$
\sn
\item[${{}}$]   $(ii) \quad$ no $\zeta \in \lambda \backslash 
\beta_{\text{Rang}(f_\alpha)}$ is $G^*$-connected to every $f_\alpha
(\gamma),\gamma \in C_\alpha$
\sn
\item[${{}}$]  $(iii) \quad \mu 
\le \beta < \alpha \Rightarrow |C_\beta \cap C_\alpha| < \kappa$.
\end{enumerate}
\mn
In stage $\alpha$ choose $B'_\alpha \subseteq B_\alpha$ of order type $\mu$
such that $(\forall \beta)[\mu \le \beta < \alpha \Rightarrow 
\sup(B'_\alpha \cap C_\beta) < \mu)$, that is 
$B'_\alpha \cap C_\beta$ is bounded in $\mu$ equivalently $C_\beta$
equivalently in $B'_\alpha$ for $\beta < \alpha$; this is possible by
diagonalization, just remember cf$(\mu) = \text{ cf}(\kappa)$ and 
$\mu > \kappa$ and clause $\circledast(i)$.  

Now there is $C$ satisfying
\mn
\begin{enumerate}
\item[$(*)^\alpha_C$]   $C \subseteq B'_\alpha$ is unbounded of
order type $\kappa$ such that no $\zeta \in
\lambda \backslash \beta_{\text{Rang}(f_\alpha)}$ is $G^*$-connected
to every $f_\alpha(\gamma)$ for $\gamma \in C$. 
\end{enumerate}
\mn
[Why?  Otherwise for every such $C$ there is a counterexample
$\gamma_C$ and we can easily choose $C_{\alpha,i}$ 
by induction on $i < \lambda$ such that:
\mn
\begin{enumerate}
\item[$\boxtimes(i)$]   $C_{\alpha,i} \subseteq B'_\alpha$
\sn
\item[$(ii)$]   sup$(C_{\alpha,i}) = \sup(B'_\alpha) = \mu$
\sn
\item[$(iii)$]  otp$(C_{\alpha,i}) = \kappa$
\sn
\item[$(iv)$]   $(\forall j < i)[\kappa > |C_{\alpha,i} \cap
C_{\alpha,j}|]$
\sn
\item[$(iv)^+$]   moreover, if $j < i$ then $\kappa > |C_{\alpha,i}
\cap \cup\{\zeta < \mu:f_\alpha(\zeta)$ is $G^*$-connected to
$\gamma_{C_{\alpha,0} \cup C_{\alpha,j}}\}|)$.
\end{enumerate}
\mn
This is easy: for clause $(iv)^+$ note that for $C = C_{\alpha,j} \cup
C_{\alpha,0}$ by the choice of $\gamma_C$ we have  
$\gamma_C \ge \beta_{\text{Rang}(f_\alpha)}$ hence by the choice of
$\beta_{\text{Rang}(f_\alpha)}$ clearly $D_C =: \{i <
\mu:f_\alpha(i)$ is well defined and $G^*$-connected to $\gamma_C\}$ has 
cardinality $< \mu$, so we can really carry the induction on $i <
\lambda$, that is any $C \subseteq B'_\alpha$ unbounded in $\mu$ of
order type $\kappa$ such that $j < i \Rightarrow |C \cap
D_{C_{\alpha,j} \cup C_{\alpha,0}}| < \kappa$ will do.  

\noindent
Let $A_0 = C_{\alpha,0},A_1 = \{\gamma_{C_{\alpha,0} \cup
C_{\alpha,1+i}}:i < \lambda\}$ they form a complete bipartite subgraph of
$G^*$ by the definition of $\gamma_{C_{\alpha,0} \cup C_{\alpha,i}}$
and $|C_{\alpha,0}| = \kappa = |A_0|$ (by $(iii)$ of $\boxtimes$) and
$|A_1| = \lambda$ (the last: by $(iv)^+$), contradiction.  So there is
$C$ such that $(*)^\alpha_C$.]

Choose $C_\alpha$ as any such $C$ such that $(*)^\alpha_C$.   Lastly define $G$

$$
U^G = \mu
$$

$$
V^G = \lambda \backslash \mu
$$

$$
R^G = \{(\beta,\alpha):\alpha \in V^G,\beta \in C_\alpha\}.
$$
\mn
Clearly $G \in {\gH}^{\text{sbp}}_{\lambda,\theta,\kappa}$ recalling
$\kappa < \theta$ and $\alpha_1 \ne \alpha_2 \Rightarrow |C_{\alpha_1}
\cap C_{\alpha_2}| < \kappa$.  Suppose
toward contradiction that $f:\lambda \rightarrow \lambda$ is a weak
embedding of $G$ into $G^*$, hence the set $Y = \{\alpha <
\lambda:\alpha \ge \mu$ and $f_\alpha = f \restriction \mu\}$ is
unbounded in $\lambda$ and \wilog \, $\alpha \in Y \Rightarrow
\beta_{\text{Rang}(f_\alpha)} = \beta^*$, i.e. is constant.  As $f$
is one to one for every $\alpha \in Y$ large enough, $f(\alpha) \in
(\beta^*,\lambda)$ and we get easy contradiction to clause
$\circledast(ii)$ for $\alpha$ and we are done.
(Note that we can add $\lambda$ nodes to $U^G$).
\end{PROOF}

For the next claim, we need another pair of definitions:
\begin{definition}
\label{2.12m}  
1) ${\gH}^*_\lambda$ is the class of $G = 
(V^G,R^G,P^G_i)_{i < \lambda}$ where $(V^G,R^G)$ is a graph, 
$|V^G| = \lambda$ and $\langle P^G_i:i < \lambda \rangle$ is a
partition of $V^G$.

\noindent
2) We say $f$ is a strong embedding of $G_1 \in {\gH}^*_\lambda$ into
$G_2 \in {\gH}^*_\lambda$ \when \,  it strongly embeds
$(V^{G_1},R^{G_1})$ into $(V^{G_2},R^{G_2})$ mapping $P^{G_1}_i$
into $P^{G_2}_i$ for $i < \lambda$.

\noindent
3) $G \in {\gH}^*_\lambda$ is ste-universal is defined naturally.
\end{definition}

\begin{definition}
\label{2.12u}
For $\kappa \le \mu,\bold U_\kappa(\mu) = \text{ min}\{|P|:P \subseteq
[\mu]^\kappa$ and $(\forall A \in [\mu]^\kappa)(\exists B \in P)(|A
\cap B| = \kappa)\}$.
\end{definition}  

\begin{remark}
\label{2.12uu}
If $\mu$ is a strong limit cardinal and cf$(\mu) < \text{ cf}(\kappa)
\le \kappa < \mu$, then $\bold U_\kappa(\mu) = \mu$.
\end{remark}

\begin{claim}
\label{2.12} 
1) Assume 
$(\lambda \ge \theta \ge \kappa \ge \aleph_0$ and)
\mn
\begin{enumerate}
\item[$(a)$]   $\lambda = \mu^+ = 2^\mu$
\sn
\item[$(b)$]  $\kappa < \mu$ and {\rm cf}$(\kappa) \ne \text{\rm cf}(\mu)$
\sn
\item[$(c)$]   $2^\kappa \le \mu$ and $\bold U_\kappa(\mu) = \mu$
\sn
\item[$(d)$]   $(i) \quad \kappa > { \text{\rm cf\/}}(\mu)$
\underline{or} 
\sn
\item[${{}}$]  $(ii) \quad \theta < \lambda$ \underline{or} 
\sn
\item[${{}}$]  $(iii) \quad  \kappa < \text{\rm cf}(\mu)$ and there 
are $C^*_\alpha \subseteq \mu$ of order type $\kappa$ 
for $\alpha < \lambda$ 

\hskip25pt such 
that $u \in [\lambda]^\lambda \Rightarrow \text{\rm otp}
[\cup\{C^*_\alpha:\alpha \in u\}] > \kappa$.
\end{enumerate}
\mn
\Then \, in ${\gH}_{\lambda,\theta,\kappa}$ there is no
we-universal even for the bipartite graphs in ${\gH}_{\lambda,\theta,\kappa}$.

\noindent
2) In part (1) if we replace clause (d) from the assumption by $(d)_1$
or $(d)_2$ where 
\mn
\begin{enumerate}
\item[$(d)_1$]  $\mu^\kappa \ge \lambda$
\sn
\item[$(d)_2$]  $\theta = \lambda$ and among the graphs or cardinal
$\mu$ there is no ste-universal
\sn
\item[$(d)_3$]  $\theta = \lambda$ and in $\cH^*_\mu$ there is no
ste-universal, \then \,still there is no 
{\rm ste}-universal, even for the bipartite graphs in
${\gH}_{\lambda,\theta,\kappa}$. 
\end{enumerate}
\mn
3) If (a) + (b) of part (1) and $2^{< \mu} = \mu = \mu^\kappa \wedge
\theta = \lambda$ \then \, there is a {\rm ste}-universal.

\noindent
4) If (a),(b) of part (1) and (c),(d) below, \then \, there is a 
{\rm ste}-universal in ${\gH}_{\lambda,\theta,\kappa}$
\mn
\begin{enumerate}
\item[$(c)$]    $\mu = \mu^\kappa$
\sn
\item[$(d)$]  in ${\gH}^*_\mu$ there is a universal, see \ref{2.12m}(1).
\end{enumerate}
\end{claim}

\begin{remark}
\label{2.12d}
Note that part (2) is not empty: if $\mu$ is strong limit singular,
$2^\mu = \lambda = \mu^+,\chi = \chi^{< \chi} < \mu$ and $\bbP$ is the
forcing of adding $\mu$ Cohen subsets to $\chi$, then in $\bold
V^{\bbP}$ clause $(d)_2$ holds.
\end{remark}

\begin{PROOF}{\ref{2.12}}
1)  Let $G^* \in {\gH}_{\lambda,\theta,\kappa}$ and \wilog
\, $V^{G^*} = \lambda$.  As in the proof of \ref{2.10} using
assumption (c) there is a club $C$ of $\lambda$ such that
\mn
\begin{enumerate}
\item[$(i)$]   $\delta \in C,\delta \le \beta < \lambda \Rightarrow
\kappa \ge \text{ otp}\{\gamma < \delta:\gamma$ is $G^*$-connected to
$\beta\}$
\sn
\item[$(ii)$]   $\delta \in C$, cf$(\delta) \ne \text{ cf}(\kappa),
\delta \le \beta < \lambda \Rightarrow \kappa > \text{ otp}\{\gamma < \beta:
\gamma$ is $G^*$-connected to $\beta\}$
\sn
\item[$(iii)$]   $\mu^2$ divides $\delta$ for every $\delta \in C$.
\end{enumerate}
\mn
Let $S =: \{\delta \in C:\text{cf}(\delta) = \text{ cf}(\kappa)\}$; as
we have $\diamondsuit_S$, see \cite{Sh:922}, so let
$\bar f = \langle f_\delta:\delta \in S \rangle,
f_\delta \in {}^\delta \delta$ be a one-to-one function 
such that $(\forall f \in {}^\lambda \lambda)(\exists^{\text{stat}}
\delta \in S)[f$ is one-to-one $\Rightarrow f_\delta = f 
\restriction \delta]$.  For $\delta \in S$ let
$\beta_\delta = \text{ Min}(C \backslash (\delta +1))$, and for 
$i \in [\delta,\beta_\delta)$ we let $a_{\delta,i} = \{\gamma < 
\delta:f_\delta(\gamma)$
is $G^*$-connected to $i\}$, so otp$(a_{\delta,i}) \le \kappa$ by the
choice of $C$, and let $B_\delta = \{i:\delta \le i <
\beta_\delta$ and $|a_{\delta,i}| \ge \kappa\}$.

Now for $\delta \in S$  we choose 
$a^*_\delta \subseteq \delta$ unbounded of order type
$\kappa$ such that $(\forall i \in B_\delta)(a^*_\delta \nsubseteq
a_{\delta,i})$.

[Why?  First assume $(d)(i)$, i.e. $\kappa > \text{ cf}(\mu)$ let
$[\delta,\beta_\delta) =
\bigcup\limits_{\xi < \text{ cf}(\mu)} A_{\delta,\xi},|A_{\delta,\xi}| < \mu,
A_{\delta,\xi}$ is $\subseteq$-increasing continuous with $\xi$ and let
$\langle \gamma_{\delta,\varepsilon}:\varepsilon < \kappa \rangle$ be
increasing continuous with limit $\delta$ satisfying
$\mu|\gamma_{\delta,\varepsilon}$
(remember $\delta \in S \Rightarrow \mu^2|\delta$).  Now choose
$\gamma^*_{\delta,\varepsilon} \in [\gamma_{\delta,\varepsilon},
\gamma_{\delta,\varepsilon +1}) \backslash \cup\{a_{\delta,i}$ : for
some $\zeta < \text{ cf}(\mu),\varepsilon = \zeta$ mod cf$(\mu)$ and
$i \in A_{\delta,\varepsilon}\}$ for $\varepsilon < \kappa$ and let
$a^*_\delta = \{\gamma^*_{\delta,\varepsilon}:\varepsilon < \kappa\}$, it
is as required.

Second assume case (ii) of clause (d) of the assumption, so $\theta <
\lambda$ hence $\theta \le \mu$.  For $\delta \in S$ we choose 
a sequence $\bar C^\delta = \langle C_{\delta,i}:i <
\mu \rangle$ of pairwise disjoint sets, 
$C_{\delta,i}$ an unbounded subset of $\delta \backslash S$ of order
type $\kappa$, always exist as $\mu^2|\delta$ (we could have asked 
moreover that $f_\delta \restriction C_{\delta,i}$ is increasing with limit
$\delta$.
Now if $f:\lambda \rightarrow \lambda$ is one-to-one 
then $\{\delta \in S:f_\delta = f \restriction
\delta$ and for $f_\delta$ we can choose $\bar C^\delta\}$ is a
stationary subset of $\lambda$ so this is O.K. but not necessary).  If
for some $i < \mu$ the set $C_{\delta,0} \cup C_{\delta,1+i}$ is as required
on $a^*_\delta$, fine, otherwise for every $i < \mu$ 
there is $\gamma_i < \lambda$ which is $G^*$-connected to
every $y \in \text{ Rang}(f \restriction (C_{\delta,0} \cup
C_{\delta,1+i}))$.  As any $\gamma_i$ is $G^*$-connected to $\le
\kappa$ ordinals $< \delta$ and $\langle C_{\delta,i}:i < \mu\rangle$
are pairwise disjoint, clearly $|\{j:\gamma_j = \gamma_i\}| 
\le \kappa$ hence we can find $Y \subseteq \mu$ such that $\langle
\gamma_i:i \in Y \rangle$ is with no repetitions and $|Y| = \theta$.
So $A_0 = \text{ Rang}(f \restriction C_{\delta,0}),
A_1 = \{\gamma_i:i \in Y\}$ exemplify that a complete 
$(\kappa,\theta)$-bipartite graph can
be weakly embedded into $G^*$, contradiction. 

Lastly, the case clause (iii) of clause (d) holds.  let $\langle
\gamma_{\delta,\varepsilon}:\varepsilon < \kappa\rangle$ be an
increasing limit with $\delta$ such that
$\mu|\gamma_{\delta,\varepsilon}$; let $\langle C^*_i:i < \lambda\rangle$
be as in clause (iii) of (d) of the assumption and let $C_{\delta,i}
:= \{\beta +1$: for some $\varepsilon < \kappa$ we have
$\gamma_{\delta,\varepsilon} \le \beta <
\gamma_{\delta,\varepsilon} + \mu$ and $\beta -
\gamma_{\delta,\varepsilon} \in C^*_i$ and otp$(C^*_i \cap (\beta -
\gamma_{\delta,\varepsilon})) < \varepsilon\}$.

Lastly, repeat the proof of ``Second...".]

Lastly define the bipartite graph $G$ by 
$V^G = \lambda,R^G = \{(\gamma,\delta):
\delta \in S,\gamma \in a^*_\delta\}$.  Easily 
$G \in {\gH}^{\text{sbp}}_{\lambda,\theta,\kappa}$ 
and is not weakly embeddable into $G^*$ by the choice of $\bar f$. 

\noindent
2) Let $G^* \in \gH_{\lambda,\theta,\kappa},V^{G^*} = \lambda$.
We choose the club $C$, the set $S$ and the sequence $\bar f = \langle
f_\delta:\delta \in S \rangle$ as in the proof of part (1).  We shall
choose $\langle a_{\delta,i}:\delta \in S,i < \mu \rangle$ and 
define $G$ by $V^G = \lambda,E^G = \bigcup\limits_{\delta \in S}
E_\delta,E_\delta = \{\{\gamma,\delta +i\}:\gamma \in a^\delta_i,i < 
i_\mu\} \cup \{\{\delta +i,\delta +j\}:(i,j) \in R_\delta \subseteq
\mu \times \mu\}$.  Naturally $a^\delta_i$ is an unbounded subset of
$\delta$ of order-type $\kappa$.

Now it is sufficient to find for $\delta \in S$ an unbounded subset $C
= C_\delta$ of $\delta$ of order type $\kappa$ such that for no $\gamma
= \gamma_C < \lambda$ do we have $(\forall \beta < \delta)(\beta \in C
\Leftrightarrow f_\delta(\beta)R^{G^*} \gamma)$, in this case
$i_\delta =1$.  If this fails then such $\gamma_C$ is well defined
for any unbounded $C \subseteq \delta$ of order type $\kappa$;
 $C_{\delta,i} \subseteq \delta$ unbounded of order type
$\kappa$, pairwise distinct for $i < \lambda$ and $C_{\delta,1+i} \cap
C_{\delta,0} = \emptyset$; then $A_0 =:
\{f_\delta(\beta):\beta \in C_{\delta,0}\},A_1 =: \{\gamma_{C_{\delta,0} \cup
C_{\delta,1+i}}:i < \lambda\}$ exemplifies that the complete
$(\kappa,\theta)$-bipartite graph can be weakly embedded into $G^*$,
contradiction.

Clearly $(d)_2 \Rightarrow (d)_3$ so \wilog \, $(d)_3$ holds.  For
$\delta \in S$ let $\langle C_{\delta,j}:j \le \mu\rangle$ be a
sequence of distinct subsets of $\mu$ which include $\{\{\alpha <
\delta:f_\delta(\alpha)R^{G^*}(\delta +i)\}:i < i_\delta\}$ and let
$M_\delta$ be the model expanding $G^* \rest [\delta,\delta +
i_\delta]$ by $P^M_j = \{\delta +i:i < i_\delta$ and $(\forall \alpha
< \delta)[f_\delta(\alpha)R^{G^*}(\delta + i) \equiv \alpha \in
C_{\delta,j}]\}$.  As $M_\delta$ is not universal in $\gH^*_\mu$, so
let $N_\delta \in \gH^*_\mu$ witness this; \wilog \, the universe of
$N$ in $[\delta,\delta + \mu)$, and let $a_{\delta,j} =
\{f_\delta(\alpha):\alpha \in C_{\delta,j}\}$ and $R_\delta =
R^{N_\delta}$.

Now check.

\noindent
3) It suffices to prove that the assumptions of part (4) holds, the
 non-trivial part is clause (d) there, i.e. $\gH^*_\lambda$ has a
 universal member.  But $2^{\le \mu} = \mu$ so either $\mu$ is regular
 so $\mu = \mu^{< \mu}$ or $\mu$ is strong limit singular and in both
 cases this holds by Jonsson or see \cite{Sh:88r}.

\noindent
4) We choose $G_\alpha$ for $\alpha \le \lambda$ by induciton on
   $\alpha$ such that
\mn
\begin{enumerate}
\item[$\boxplus$]  $(a) \quad G_\alpha$ is a graph with set of nodes
$(1+\alpha)\mu$
\sn
\item[${{}}$]  $(b) \quad$ if $\beta < \alpha$ then $G_\beta$ is an
induced subgraph of $G$
\sn
\item[${{}}$]  $(c) \quad$ if $\alpha = \beta +1$ and $G$ is a graph
with $\mu$ nodes and id$_{G_\beta}$ is a strong 

\hskip25pt embedding of $G_\beta$ into $G$ such that $x \in V^G \backslash
V^{G_\beta} \Rightarrow$ 

\hskip25pt ($x$ is connected to $\le \kappa$ nodes of
$G_\beta$) \then \, there is a 

\hskip25pt strong embedding of $G$ 
into $G_\alpha$ which extends id$_{G_\beta}$.
\end{enumerate}
\mn
The construction is possible by clause (d) of the assumption.  Now as
in the proof of \ref{2.11} $G_\lambda \in
\gH_{\lambda,\theta,\kappa}$ is ste-universal.
\end{PROOF}

\begin{remark}
\label{2.12A}  
1) In the choice of $\bar f$ (in the proof of \ref{2.12}) 
we can require that for every $f \in {}^\lambda \lambda$ the set
$\{\delta \in S:f_\delta = f \restriction \delta$ and $\delta \cap
\text{ Rang}(f) = \text{ Rang}(f_\delta)\}$ is stationary and so 
deal with copies of
the complete $(\kappa,\theta)$-bipartite graph with the $\theta$ part
after the $\kappa$ part. 

\noindent
2) Probably we can somewhat weaken assumption (c).
\end{remark}

\begin{theorem}
\label{2.14}  
Assume $\lambda \ge \theta \ge \kappa \ge
\aleph_0$ and $\lambda = 2^\mu = \mu^+$ and $2^{<\mu} = \mu$. 

\noindent
1) In ${\gH}_{\lambda,\theta,\kappa}$ there is a we-universal member \Iff \,
$\mu^\kappa = \mu \wedge \theta = \lambda$ \Iff \,  there is no $G^* \in 
{\gH}_{\lambda,\theta,\kappa}$ we-universal for $\{G^{[\text{\rm gr}]}:
G \in {\gH}^{\text{\rm sbp}}_{\lambda,\theta,\kappa}\}$. 

\noindent
2) In ${\gH}_{\lambda,\theta,\kappa}$ there is ste-universal,
\Iff \, $\mu^\kappa = \mu \wedge \theta = \lambda$ \underline{iff}
there is no $G^* \in {\gH}_{\lambda,\theta,\kappa}$
ste-universal for $\{G^{[\text{\rm gr}]}:G \in
{\gH}^{\text{\rm sbp}}_{\lambda,\theta,\kappa}\}$.
\end{theorem}

\begin{PROOF}{\ref{2.14}}
1) The second iff we ignore as in each case the same
claims cited give it too \underline{or} use \ref{2.16} below.  We first
prove that there is a we-universal except possibly when $\mu^\kappa =
\mu \wedge \theta = \lambda$.

Proving this claim, whenever we point out a case is
resolved we assume that it does not occur.  We avoid using $2^{< \mu}
= \mu$ when we can.  

If $\lambda = \lambda^{\theta^+}$ then by \ref{2.5}(3) we have
Ps$'(\lambda,\theta^+ \times \kappa)$ so by \ref{2.4} +
\ref{2.5B}(2) there is no we-universal; hence we can assume that
$\lambda < \lambda^{\theta^+}$ so (as $\lambda = \lambda^{< \lambda}$)
clearly $\lambda \le \theta^+$ hence $\lambda = \theta \vee \lambda =
\theta^+$ that is $\theta = \lambda \vee \theta = \mu$.

If $\kappa = \theta$ then by \ref{2.8}(2) there is no we-universal,
so we can assume that $\kappa \ne \theta$ hence $\kappa < \theta \le
\lambda$, hence $\lambda = \lambda^\kappa$ so by \ref{2.5}(3) we
have Ps$'(\lambda,\kappa)$ hence by \ref{2.5}(7) we have
Ps$(\lambda,\kappa)$.  So if $\theta
= \kappa^+$ then by \ref{1.3} more exactly, \ref{2.5B}(1)
there is no we-universal so without loss of generality
$\kappa^+ < \theta$ hence $\kappa < \mu$.
If cf$(\kappa) = \text{ cf}(\mu)$ then by \ref{2.11} we are 
done except if 
\mn
\begin{enumerate}
\item[$(*)_1$]   clause (c) of \ref{2.11} fails, (and
cf$(\kappa) = \text{ cf}(\mu),\kappa < \mu$)

is impossible as $2^{< \mu} = \mu$.
\end{enumerate}
\mn
But (c) of \ref{2.11} so we can assume
\mn
\begin{enumerate}
\item[$(*)_2$]   cf$(\kappa) \ne \text{ cf}(\mu)$,
\end{enumerate}
\mn
so as $2^{< \mu} = \mu,\kappa < \mu$ we get
\mn
\begin{enumerate}
\item[$(*)_3$]   $\bold U_\kappa(\mu) = \mu$ and $2^\kappa \le \mu$.
\end{enumerate}
\mn
Now we try to apply \ref{2.12}(1), so we can assume that we cannot; but
clauses (a)-(c) there hold hence clause (d) there fails.
So $\kappa \le \text{\rm cf}(\mu) \wedge \theta \ge \lambda$
 (recalling sub-clauses (i),(ii) of \ref{2.12}(1)(d)) as
cf$(\kappa) \ne \text{ cf}(\mu)$ by $(*)_2$ and $\theta \le \lambda$ 
we have $\kappa < \text{ cf}(\mu)$ and $\theta = \lambda$.  As $2^{<
\mu} = \mu$ this implies 
$\mu^\kappa = \mu$ and $\theta = \lambda = \mu^+$ as promised.  All
this gives the implication $\Rightarrow$; the other direction by
\ref{2.10} gives there is a we-universal. 

\noindent
2) By part (1) and \ref{0.5}(3), the only open case is $\mu^\kappa
= \mu$ and $\theta = \lambda = \mu^+$ then Claim \ref{2.12}(3),(4) 
applies (clause (c) there follows from $\mu = \mu^\kappa$).
\end{PROOF}

Recalling Definition \ref{0.4}
\begin{claim}
\label{2.15}  
The results in \ref{2.14} hold for
${\gH}^{\text{\rm sbp}}_{\lambda,\theta,\kappa}$ and for 
${\gH}^{\text{\rm bp}}_{\lambda,\theta,\kappa}$.
\end{claim}

\begin{PROOF}{\ref{2.15}}
The "no universal'' clearly holds by \ref{2.14}, so
we need the ``positive results", and we are done by \ref{2.16} below.
\end{PROOF}

\begin{claim}
\label{2.16}   The results of \ref{2.6}, \ref{2.7},
\ref{2.10} and \ref{2.12}(3),(4) hold for 
${\gH}^{\text{\rm sbp}}_{\lambda,\theta,\kappa}$ and
${\gH}^{\text{\rm bp}}_{\lambda,\theta,\kappa}$.
\end{claim}

\begin{PROOF}{\ref{2.16}}
  In all the cases the isomorphism and embeddings preserve
$``x \in U^G",``y \in V^G"$.

For ${\gH}^{\text{sbp}}_{\lambda,\theta,\kappa}$, in \ref{2.10} we
redefine $G^*$ as a bipartite graph (recalling $\langle a^\delta_i:i <
\mu \rangle$ lists $\{a \subseteq \delta:\text{ otp}(a) \le \kappa$ and
if equality holds then $\delta = \sup(a)\}$ for 
$\delta < \lambda$ divisible by $\mu$)

$$
U^{G^*} = \{2 \alpha:\alpha < \lambda\}
$$

$$
V^{G^*} = \{2 \alpha +1:\alpha < \lambda\}
$$

\begin{equation*}
\begin{array}{clcr}
R^{G^*} = &\{(2 \alpha,2 \beta +1):\text{for some } \delta < \lambda
\text{ divisible by } \mu \text{ we have } 2\alpha,2 \beta +1 \in
[\delta,\delta + \mu]\} \\
  &\cup\{(2 \alpha,\delta + 2i+1):
\delta < \lambda \text{ divisible by } \mu,i < \mu,2 \alpha <
\delta,\alpha \in a^\delta_i\} \\
  &\cup\{(\delta +2i,2 \beta + 1):\delta < \lambda \text{ divisible by }
\mu,2 \beta + 1 < \delta,i < \mu,\beta \in a^\delta_i\}
\end{array}
\end{equation*}
\mn
The proof is similar.
\sn
For ${\gH}^{\text{bp}}_{\lambda,\theta,\kappa}$, \ref{2.10} we
redefine $G^*$ 

$$
U^{G^*} = \{2 \alpha:\alpha < \lambda\}
$$

$$
V^{G^*} = \{2 \alpha +1:\alpha < \lambda\}
$$

\begin{equation*}
\begin{array}{clcr}
R^{G^*} = &\{(2 \alpha,2 \beta +1):\text{for some } \delta < \lambda
\text{ divisible by } \mu,\{2 \alpha,2 \beta +1\} \subseteq
[\delta,\delta + \mu]\} \\
  &\cup\{(2 \alpha,\delta + 2i+1):
\delta < \lambda \text{ divisible by } \mu,i < \mu,2 \alpha <
\delta \text{ and } \alpha \in a^\delta_i\} \\
  &\cup\{(2 \alpha +1,2 \beta):2 \alpha +1 < 2 \beta\}.
\end{array}
\end{equation*}
\mn
The proof of \ref{2.6}, \ref{2.7} for 
${\gH}^{\text{sbp}}_{\lambda,\theta,\kappa}$ is similar 
to that of \ref{2.6}, \ref{2.7}.  The $G_\eta$ is from
${\gH}^{\text{sbp}}_{\lambda,\ell g(\eta)}$ so the isomorphism preserve
the $x \in U^G,y \in V^G$.  For ${\gH}^{\text{bp}}_{\lambda,
\theta,\kappa}$ without loss of generality $\kappa \ne \theta$ hence
$\kappa < \theta$ (otherwise this falls under the previous case).  We
repeat the proof of the previous case carefully; making the following
changes, say for \ref{2.6}, $\langle V^G_i:i < \sigma \rangle$ is
increasing continuous with union $V^G,\langle U^G_i:i < \sigma
\rangle$ increasing continuous with union $U^G$.
\mn
\begin{enumerate}
\item[$(*)'_1$]   if $x \in V^G \backslash V^G_{i+1}$ then $\kappa >
|\{y \in U^G_i:y$ is $G$-connected to $x\}|$.  
\end{enumerate}
\mn
We leave \ref{2.12}(3),(4) to the reader.  
\end{PROOF}

\section {More accurate properties} 

\begin{definition}
\label{4.1}   Let Q$(\lambda,\mu,\sigma,\kappa)$ 
mean: there are $A_i \in [\mu]^\sigma$ for 
$i < \lambda$ such that for every $B \in
[\mu]^\kappa$ there are $< \lambda$ ordinals $i$ such that $B
\subseteq A_i$. 
\end{definition}

\begin{definition}
\label{4.2}  
1) For $\mu \ge \kappa$ let set$(\mu,\kappa) 
= \{A:A$ is a subset of $\mu \times \kappa$ of
cardinality $\kappa$ such that $i < \kappa \Rightarrow \kappa > |\{A
\cap (\mu \times i)\}|$ and let set$(\mu,\kappa) = [\kappa]^\kappa$
for $\mu < \kappa$. 

\noindent
2) Assume $\lambda \ge \theta \ge \kappa,\lambda \ge \mu$.  
Let Qr$_{\text{\rm w}}(\lambda,\mu,\theta,\kappa)$ 
mean that some $\bar A$ exemplifies it, which means
\mn
\begin{enumerate}
\item[$(a)$]   $\bar A = \langle A_i:i < \alpha \rangle$
\sn
\item[$(b)$]   $A_i \in \text{ set}(\mu,\kappa)$ for $i < \alpha$
\sn
\item[$(c)$]  $\bar A$ is $(\kappa,\theta)$-free which means 
$(\forall A \in \text{ set}(\mu,\kappa))(\exists^{< \theta} 
i < \alpha)(A \subseteq A_i)$
\sn
\item[$(d)$]  $\alpha \le \lambda$
\sn
\item[$(e)$]   if $\bar A' = \langle A'_i:i < \alpha' \rangle$
satisfies clauses (a),(b),(c),(d), then for
some one to one function $\pi$ from $\bigcup\limits_{i < \alpha'} A'_i$
into $\bigcup\limits_{i < \alpha} A_i$ and
one-to-one function $\varkappa$ from $\alpha'$ to $\alpha$ 
(or $\kappa$ to $\kappa$) we have $i < \alpha' \Rightarrow 
\pi(A'_i) \subseteq A_{\varkappa(i)}$.
\end{enumerate}
\mn
3) Qr$_{\text{\rm st}}(\lambda,\mu,\theta,\kappa)$ is defined similarly except
that we change clause (e) to (e)$^+$ demanding $\pi(A'_i) = A_{\varkappa(i)}$. 
Let Qr$_{\text{pr}}(\lambda,\mu,\theta,\kappa)$ be defined similarly
omitting clause (e). 

\noindent
4) Assume $\lambda \ge \theta \ge \kappa,\lambda \ge \mu$ and $x \in
\{w,st\}$.  Let NQr$_x(\lambda,\mu,\theta,\kappa)$ mean that 
Qr$_x(\lambda,\mu,\theta,\kappa)$ fails.
\end{definition}

\begin{claim}
\label{4.3}  
1) Assume {\rm NQr}$_w(\lambda,\mu,\theta,\kappa)$ 
and $\lambda = \lambda^{\mu + \kappa}$ and $(\lambda > \mu^\kappa) 
\vee (\mu < \kappa)$.  \Then
\mn
\begin{enumerate}
\item[$(a)$]   in ${\gH}_{\lambda,\theta,\kappa}$ there is no
{\rm we}-universal member
\sn
\item[$(b)$]  moreover, for every $G^* \in 
{\gH}_{\lambda,\theta,\kappa}$ there is a member of 
${\gH}^{\text{\rm sbp}}_{\lambda,\theta,\kappa}$ not weakly embeddable into it.
\end{enumerate}
\mn
2) Assume {\rm NQr}$_{\text{\rm st}}(\lambda,\mu,\theta,\kappa)$ 
and $\lambda = \lambda^{\mu + \kappa}$ and 
$(\lambda > \mu^\kappa) \vee (\mu \le \kappa)$.  \Then
\mn
\begin{enumerate}
\item[$(a)$]   in ${\gH}_{\lambda,\theta,\kappa}$ there is no
{\rm ste}-universal member
\sn
\item[$(b)$]   moreover, for every $G^* \in {\gH}_{\lambda,\theta,\kappa}$ 
there is a member of ${\gH}^{\text{\rm sbp}}_{\lambda,\theta,\kappa}$ not
strongly embeddable into it. 
\end{enumerate}
\mn
3) In parts (1), (2) we can weaken the assumption $\lambda =
\lambda^{\mu + \kappa}$ to
\mn
\begin{enumerate}
\item[$\otimes$]   $\lambda = \lambda^\kappa > \mu$ and there is
${\cF} \subseteq \{f:f$ a partial one to one function from
$\lambda$ to $\lambda,|{\text{\rm Dom\/}}(f)| = \mu\}$ of cardinality
$\lambda$ such that for every $f^* \in {}^\lambda \lambda$ there
\footnote{we can add ``there are
$\lambda$ functions $f \in {\cF},f \subseteq f^*$, with pairwise
disjoint domains", and possibly increasing ${\cF}$ we get it}  is
$f \in {\cF},f \subseteq f^*$.
\end{enumerate}
\end{claim}

\begin{PROOF}{\ref{4.3}}
1), 2)  Let $x = w$ for part (1) and $x =$ st for part (2).

Now suppose that $G^* \in {\gH}_{\lambda,\theta,\kappa}$ and
\wilog \, $V^{G^*} = \lambda$ and we shall construct $G \in 
{\gH}^{\text{sbp}}_{\lambda,\theta,\kappa}$ not $x$-embeddable into
$G^*$ (so part (b) will be proved, and part (a) follows).
\bigskip

\noindent
\underline{Case 1}:  $\mu < \kappa$ so set$(\mu,\kappa) = [\kappa]^\kappa$.

Similar to the proof of \ref{1.4}.  Let $\bar f = \langle
f_\eta:\eta \in {}^\kappa \lambda \rangle$ be a simple black box for
one to one functions.  It means that each $f_\eta$ is a one to one
function from $\{\eta \rest i:i < \kappa\}$ into $\lambda$, such that
for every $f:{}^{\kappa >} \lambda \rightarrow \lambda$ for some $\eta
\in {}^\kappa \lambda$ we have $f_\eta \subseteq f$.
We define the bipartite graph $G$ as follows:
\mn
\begin{enumerate}
\item[$(*)$]   $(a) \quad U^G = {}^{\kappa >}\lambda$ and $V^G = 
({}^{\kappa}\lambda) \times \mu$
\sn
\item[${{}}$]  $(b) \quad R^G = \cup\{R^G_\eta:\eta \in {}^\kappa \lambda\}$
where $R^G_\eta \subseteq \{(\eta \restriction
\varepsilon,(\eta,i)):\varepsilon < \kappa,i < \mu\}$.
\end{enumerate}
\mn
Now for each $\eta \in {}^\kappa \lambda$, we choose
$\langle \beta_{\eta,i}:i < \alpha_\eta \rangle$ listing without
repetitions the set $\{\beta < \lambda:\beta$ is $G^*$-connected to $\kappa$
members of Rang$(f_\eta)\}$ and \wilog \, $\alpha_\eta =
|\alpha_\eta|$ and $A_{\eta,i} = \{\varepsilon <
\kappa:\beta_{\eta,i},f_\eta(\eta \restriction \varepsilon)$ are
$G^*$-connected$\}$.

As $G^* \in \gH_{\lambda,\theta,\kappa}$ 
clearly $A \in [\kappa]^\kappa \Rightarrow |\{i < \alpha_\eta:A
\subseteq A_{\eta,i}\}| < \theta$ hence $\alpha_\eta = |\alpha_\eta|
\le 2^\kappa + \theta$ but (see Definition \ref{4.2}(2))
we have assumed $\theta \le \lambda$ and (in \ref{4.3})
$\lambda = \lambda^\kappa$ so $2^\kappa \le \lambda$ hence
$\alpha_\eta \le \lambda$; next let $\bar A_\eta = \langle
A_{\eta,i}:i < \alpha_\eta \rangle$ and as $\bar A_\eta$ cannot be a
witness for Qr$_x(\lambda,\mu,\theta,\kappa)$ but clauses (a), (b),
(c), (d) of Definition \ref{4.2}(1) hold, hence clause (e) fails so there
is $\bar A'_\eta = \langle A'_{\eta,i}:i < \alpha'_\eta \rangle$
exemplifies the failure of clause (e) of Definition \ref{4.2} with
$\bar A_\eta,\bar A'_\eta$ here standing for $\bar A,\bar A'$ there.

Let 

\[
R^G_\eta = \{(\eta \restriction
\varepsilon,(\eta,i)):\varepsilon < \kappa,i < \alpha'_\eta \text{ and }
\varepsilon \in A'_{\eta,i}\}.
\]

\mn
The proof that $G$ cannot be $x$-embedded into $G^*$ is as in the
proof of \ref{1.4}.
\bigskip

\noindent
\underline{Case 2}:  $\mu^\kappa < \lambda$ (and $\lambda = \lambda^\mu,\kappa
\le \mu < \lambda$).

First note that by the assumptions of the case
\mn
\begin{enumerate}
\item[$\boxplus_1$]   there is $\bar f = \langle f_\eta:\eta \in
{}^\kappa \lambda \rangle$ such that
\begin{enumerate}
\item[$(a)$]   $f_\eta$ is a function from $\bigcup\limits_{\varepsilon <
\kappa} (\{\eta \restriction \varepsilon\} \times \mu)$ into $\lambda$
\sn
\item[$(b)$]   if $f$ is a function from $({}^{\kappa >} \lambda)
\times \mu$ to $\lambda$ then we can find $\langle \nu_\rho:\rho \in
{}^{\kappa \ge} \lambda \rangle$ such that
\sn
\item[${{}}$]  $\qquad (i) \quad \nu_\rho \in {}^{\ell g(\rho)} \lambda$
\sn
\item[${{}}$]  $\qquad (ii) \quad \rho_1 \triangleleft \rho_2 \Rightarrow
\nu_{\rho_1} \triangleleft \nu_{\rho_2}$
\sn
\item[${{}}$]  $\qquad (iii) \quad$ if $\alpha < \beta < \lambda$ and $\rho \in
{}^{\kappa >} \lambda$ then $\nu_{\rho \char 94 \langle \alpha
\rangle} \ne \nu_{\rho \char 94 \langle \beta \rangle}$
\sn
\item[${{}}$]  $\qquad (iv) \quad f_{\nu_\rho} \subseteq f$ for $\rho \in
{}^\kappa \lambda$.
\end{enumerate}
\end{enumerate}
\mn
We commit ourselves to
\mn
\begin{enumerate}
\item[$\boxplus_2$]   $(a) \quad U^G = 
({}^{\kappa >}\lambda) \times \mu$ and $V^G =
\{(\eta,i):\eta \in {}^\kappa \lambda,i < \lambda\}$
\sn
\item[${{}}$]  $(b) \quad R^G = \cup\{R^G_\eta:\eta \in {}^\kappa
\lambda\}$ where
\sn
\item[${{}}$]  $(c) \quad R^G_\eta \subseteq \{((\eta \restriction
\varepsilon,j),(\eta,i)):j < \mu,i < \lambda,\varepsilon < \kappa\}$.
\end{enumerate}
\mn
We say $\eta \in {}^\kappa \lambda$ is $G^*$-reasonable if $f_\eta$ is
one to one and for every
$\zeta < \kappa$ and $y \in V^{G^*}$ the set $\{(\eta \restriction
\varepsilon,j):\varepsilon < \zeta,j < \mu$ and $f_\eta((\eta
\restriction \varepsilon,j))$ is $G^*$-connected to $y\}$ has cardinality
$< \kappa$.  We decide
\mn
\begin{enumerate}
\item[$\boxplus_3$]   $(a) \quad$ if $\eta$ is 
not $G^*$-reasonable then $R^G_\eta = \emptyset$
\sn
\item[${{}}$]  $(b) \quad$ if $\eta$ is $G^*$-reasonable let $\langle
\beta_{\eta,i}:i < \alpha_\eta \rangle$ list without repetitions the

\hskip25pt set $\{\beta < \lambda:\beta$ is 
$G^*$-connected to at least $\kappa$ members of 

\hskip25pt Rang$(f_\eta)\}$; and let $A_{\eta,i} 
= \{(\varepsilon,j):\varepsilon
< \kappa,j < \mu$ and $f_\eta((\eta \restriction \varepsilon,j))$

\hskip25pt is $G^*$-connected to $\beta_{\eta,i}\}$; 
clearly $A_{\eta,i} \in \text{ set}(\mu,\kappa)$ and 

\hskip25pt let $\bar A_\eta = \langle A_{\eta,i}:i < \alpha_\eta \rangle$
\sn
\item[${{}}$]  $(c) \quad$ as $\bar A_\eta$ cannot guarantee
{\rm Qr}$_x(\lambda,\mu,\theta,\kappa)$ necessarily there is 

\hskip25pt $\bar A'_\eta =
\langle A'_{\eta,i}:i < \alpha'_\eta \rangle$ exemplifying this so
$\alpha'_\eta \le \lambda$ and let 

\hskip25pt $R^G_\eta = \{((\eta \restriction
\varepsilon,j),(\eta,i)):i < \alpha'_\eta,$ and $(\varepsilon,j) \in
A'_{\eta,i}\}$.
\end{enumerate}
\mn
The rest should be clear; for every $f:{}^{\kappa >}\lambda
\rightarrow \lambda$ letting $\langle \nu_\rho:\rho \in {}^{\kappa
\ge}\lambda \rangle$ be as in $\boxtimes$ above, for some $\rho \in
{}^\kappa \lambda,\nu_\rho$ is $G^*$-reasonable.
\bigskip

\noindent
\underline{Case 3}:  $\mu = \kappa$.

Left to the reader (as after Case 1,2 it should be clear). 

\noindent
3) As in the proof of \ref{2.5}(4), it follows that there 
is $\bar f$ as needed.
\end{PROOF}  

\begin{claim}
\label{4.4}  
1) {\rm NQr}$_{\text{\rm st}}(2^\kappa,\kappa,2^\kappa,\kappa)$. 

\noindent
2) If $\lambda = \theta = 2^\kappa$ \then \, in 
${\gH}_{\lambda,\theta,\kappa}$ there is no ste-universal even for
members of ${\gH}^{\text{\rm sbp}}_{\lambda,\theta,\kappa}$.
\end{claim}
 
\begin{PROOF}{\ref{4.4}}  
1)  Think. 

\noindent
2) By part (1) and \ref{4.3}.
\end{PROOF}

\begin{claim}
\label{4.5}  
1) Assume $\kappa < \lambda$ and {\rm Qr}$_x(\lambda,1,\lambda,\kappa)$ and 
${\frak H}^{\text{\rm sbp}}_{(\lambda,\kappa),\lambda,
\kappa} \ne \emptyset$, \then \, 
${\frak H}^{\text{\rm sbp}}_{(\lambda,\kappa),\lambda,\kappa} = 
{\frak H}^{\text{\rm bp}}_{(\lambda,\kappa),\lambda,\kappa}$ has a
$x$-universal member.
\end{claim}

\begin{PROOF}{\ref{4.5}}
Read the definitions.
\end{PROOF}

\section {Independence results on existence of large almost disjoint
families}

This deals with a question of Shafir
\begin{definition}
\label{3.1} 
1) Let Pr$_2(\mu,\kappa,\theta,\sigma)$
mean: there is ${\cA} \subseteq [\kappa]^\kappa$ such that
\mn
\begin{enumerate}
\item[$(a)$]   $|{\cA}| = \mu$
\sn
\item[$(b)$]   if $A_i \in {\cA}$ for $i < \theta$ and $i \ne j
\Rightarrow A_i \ne A_j$ then $|\bigcap\limits_{i < \theta} A_i| < \sigma$.
\end{enumerate}
\mn
If we omit $\sigma$ we mean $\kappa$, if we omit $\mu$ we mean $2^\kappa$.
\end{definition}

\begin{claim}
\label{3.3}  
{\rm Pr}$_2(-,-,-,-)$ has obvious monotonicity properties.
\end{claim}

\begin{claim}
\label{3.2} 
Assume
\mn
\begin{enumerate}
\item[$(*)$]   $\sigma = \sigma^{< \sigma} < \kappa = \kappa^\sigma =
\text{\rm cf}(\kappa)$ and 
$(\forall \alpha < \kappa)(|\alpha|^{< \sigma} < \kappa)$
and $2^\kappa = \kappa^+ < \chi$ (so $\kappa^{++} \le \chi$).
\end{enumerate}
\mn
\Then \, for some forcing notion $\Bbb P$ 
\mn
\begin{enumerate}
\item[$(a)$]  $|\Bbb P| = \chi^\kappa$
\sn
\item[$(b)$]  $\Bbb P$ satisfies the $\kappa^{++}$-c.c.
\sn
\item[$(c)$]  $\Bbb P$ is $\sigma$-complete
\sn
\item[$(d)$]   $\Bbb P$ neither collapses cardinals nor changes cofinalities
\sn
\item[$(e)$]  in $\bold V^{\Bbb P}$ we have 
$2^\sigma = \chi^\sigma,2^\kappa = \chi^\kappa$
\sn
\item[$(f)$]   in $\bold V^{\Bbb P}$ we have {\rm Pr}$_2(\chi,
\kappa,\sigma^+,\sigma)$ but $\neg\text{\rm Pr}_2(\kappa^{++},\kappa,
\kappa,\theta)$ for $\theta < \sigma$ recalling \ref{3.3}. 
\end{enumerate}
\end{claim}

\begin{PROOF}{\ref{3.2}}
The forcing is as in a special case of the $\bbQ$ one in 
\cite[\S2]{Sh:918}, see history there.  Let $E$ be the following 
equivalence relation on $\chi$

$$
\alpha E \beta \Rightarrow \alpha + \kappa = \beta + \kappa.
$$
\mn
We define the partial order $\Bbb P = (P,\le)$ by

\begin{equation*}
\begin{array}{clcr}
P = \{f:&f \text{ is a partial function from } \chi \text{ to } \{0,1\} \\
  &\text{ with domain of cardinality } \le \kappa \text{ such that} \\ 
  &(\forall \alpha < \chi)(|\text{Dom}(f) \cap (\alpha/E)| < \sigma)\}
\end{array}
\end{equation*}

\begin{equation*}
\begin{array}{clcr}
f_1 \le f_2 \text{ \underline{iff} } &f_1,f_2 \in P,
f_1 \subseteq f_2 \text{ and} \\
  &\sigma > |\{\alpha \in \text{ Dom}(f_1):f_1 \restriction (\alpha/E) \ne
f_2 \restriction (\alpha/E)\}|.
\end{array}
\end{equation*}
\mn
We define two additional partial orders on $P$:

\begin{equation*}
\begin{array}{clcr}
f_1 \le_{\text{pr}} f_2 \text{ \underline{iff} } &f_1 \subseteq f_2 
\text{ and } f_1,f_2 \in \Bbb P \text{ and} \\
  &(\forall \alpha \in \text{ Dom}(f_1))[f_1 \restriction (\alpha/E) =
f_2 \restriction (\alpha/E)].
\end{array}
\end{equation*}

$$
f_1 \le_{\text{apr}} f_2 \text{ iff } f_1,f_2 \in \Bbb P,f_1 \le f_2
\text{ and Dom}(f_2) \subseteq \cup \{\alpha/E:\alpha \in \text{ Dom}(f_1)\}.
$$

We know (see there)
\mn
\begin{enumerate}
\item[$(*)_0$]   $\Bbb P$ is $\kappa^{++}$-c.c., $|\Bbb P| = \chi^\kappa$
\sn
\item[$(*)_1$]  $\Bbb P$ is $\sigma$-complete and
$(P,\le_{\text{pr}})$ is $\kappa^+$-complete, in both cases the union
of an increasing sequence forms an upper bound
\sn
\item[$(*)_2$]   for each $p,\Bbb P \restriction 
\{q:p \le_{\text{apr}} q\}$ is $\sigma^+$-c.c. of 
cardinality $\kappa^\sigma = \kappa$ and $\sigma$-complete
\sn
\item[$(*)_3$]   if $p \le r$ then for some $q,q'$ we have 
$p \le_{\text{apr}} q \le_{\text{pr}} r$ and $p \le_{\text{pr}} 
q' \le_{\text{apr}} r$; moreover $q$ is unique we denote it by
inter$(p,r)$
\sn
\item[$(*)_4$]   if $p \Vdash ``\name \tau \in
{}^\kappa\text{Ord}$" then for some $q$ we have
\begin{enumerate}
\item[$(a)$]   $p \le_{\text{pr}} q$
\sn
\item[$(b)$]   if $\alpha < \kappa$ and $q \le r$ and $r \Vdash_{\Bbb P}
``\name \tau(\alpha) = \beta"$ then inter$(q,r)
\Vdash_{\Bbb P} ``\name \tau(\alpha) = \beta"$.
\end{enumerate}
\end{enumerate}
\mn
This gives that clauses (a),(b),(c),(d) of the conclusion hold.  As for clause
(e), $2^\kappa \le \chi^\kappa$ follows from $(*)_5 + 
|\Bbb P| = \chi^\kappa$ and $2^\sigma \le \chi^\sigma$, too.  

Define:
\mn
\begin{enumerate}
\item[$(*)_5$]   $(a) \quad \name f := \cup\{p:p \in {\name
G_{\bbP}}\}$
\sn
\item[${{}}$]  $(b) \quad \Vdash_{\Bbb P} ``\name f = \cup
\name G_{\Bbb P}$ is a function from $\chi$ to $\{0,1\}"$
\sn
\item[${{}}$]  $(c) \quad$ for $\alpha < \chi$ let
$\name A_\alpha = \{\gamma < \kappa:\name f(\kappa \alpha + \gamma)=1\}$.
\end{enumerate}

Also easily
\mn
\begin{enumerate}
\item[$(*)_6$]   if $p \Vdash ``\name \tau \subseteq
\sigma$" then for some $u \in [\chi]^{\le \sigma}$ and $\sigma$-Borel
funtion $\bold B:{}^u 2 \rightarrow {}^\sigma 2$ and $q$ we have 

$p \le_{\text{pr}} q$

$q \Vdash ``\name \tau = \bold B(\name f \restriction u)"$ 
\sn
\item[$(*)_7$]   $\Vdash ``\name A_\alpha \subseteq
\kappa$ moreover $\gamma < \kappa \Rightarrow 
[\gamma,\gamma + \sigma) \cap \name A_\alpha \ne
\emptyset$ and $[\gamma,\gamma + \sigma)
\nsubseteq \name A_\alpha$ and $\alpha \ne \beta
\Rightarrow \name A_\alpha \ne \name A_\beta"$.
\end{enumerate}
\mn
[Why?  By density argument.]
\mn
\begin{enumerate}
\item[$(*)_8$]  $\cA := \{\name A_\alpha:\alpha < \chi\}$ 
exemplifies Pr$_2(\chi,\kappa,\sigma^+,\sigma)$. 
\end{enumerate}
\mn
Why?  By $(*)_7 + (*)_5(c)$, ${\cA} \subseteq [\kappa]^\kappa,|{\cA}| =
\chi$ so we are left with proving clause (b) of \ref{3.1}; its
proof will take awhile.  So toward contradiction assume that
for some $p \in \Bbb P$ and $\langle \name \beta_\zeta:\zeta 
< \sigma^+ \rangle$ we have 
\mn
\begin{enumerate}
\item[$\boxplus_1$]   $p \Vdash_{\bbP} ``{\name \beta_\zeta} < \chi,
\name \beta_\zeta \ne \name \beta_\xi$ for $\zeta < \xi < \kappa^+$
and $|\bigcap_{\zeta < \sigma^+} \, \name A_{\beta_\zeta}| \ge \sigma"$.
\end{enumerate}
\mn
By induction on $\zeta \le \sigma^+$ we choose $p_\zeta \in \bbP$ such that:
\mn
\begin{enumerate}
\item[$\boxplus_2$]   $(\alpha) \quad p_0 = p$
\sn
\item[${{}}$]  $(\beta) \quad p_\zeta$ 
is $\le_{\text{pr}}$-increasing continuous
\sn
\item[${{}}$]  $(\gamma) \quad$ there is $r'_\zeta$ such that
$p_{\zeta +1} \le_{\text{apr}} r'_\zeta$ and 
$r'_\zeta \Vdash \name \beta_\zeta = \beta^*_\zeta$
\sn
\item[${{}}$]  $(\delta) \quad$ Dom$(p_{\zeta +2}) \cap (\beta^*_\zeta/E)
 \ne \emptyset$.
\end{enumerate}
\mn
No problem because $(P,\le_{\text{pr}})$ is $\kappa^+$-complete and
$\sigma^+ < \kappa^+$ and $(*)_3$.

Let\footnote{if we define $\bbP$ such that it is only
$\kappa$-complete, we first choose $y_2,u_3,u^*,v_\zeta,v^*$ and then
$q=p_{\zeta(*)}$ for $\zeta(*) = \text{ min}\{\zeta < \sigma^+:|\zeta
\cap y| = \sigma\}$}  $q = p_{\sigma^+}$. 

We can find $r_\zeta,\beta^*_\zeta$ for $\zeta < \sigma^+$ such that
$q \le_{\text{apr}} r_\zeta$ and $r_\zeta \Vdash
\name \beta_\zeta = \beta^*_\zeta$.  Let $u_\zeta =
\text{ Dom}(r_\zeta) \backslash \text{ Dom}(q)$, so $|u_\zeta| < \sigma$ by
the definition of $\le_{\text{apr}}$. By the $\Delta$-system lemma,
recalling $\sigma = \sigma^{< \sigma}$ there is
$Y \subseteq \sigma^+,|Y| = \sigma^+$ and $u^*$ such that for $\zeta < \xi$
from $Y$ we have 
$u_\zeta \cap u_\xi = u^*$.  Without loss of generality $\zeta \in Y
\Rightarrow r_\zeta \restriction u^* = r^*$.  Let $v_\zeta = \{\gamma <
\kappa:\kappa \beta^*_\zeta + \gamma \in \text{ Dom}(r_{\zeta +2})\}$ so
$v_\zeta \in [\kappa]^{< \sigma}$.

Possibly further shrinking $Y$ \wilog \,

$$
\zeta \ne \xi \text{ from } Y \Rightarrow v_\zeta \cap v_\xi = v^*.
$$
\mn
So $v^* \in [\kappa]^{< \sigma}$ (in fact follows).

Let $\zeta(*) = \text{ Min}(Y)$. 

We claim

\[  
r_{\zeta(*)} \Vdash ``\bigcap\limits_{\zeta < \sigma^+} 
{\name A_{\name \beta_\zeta}} \subseteq v^*".
\]
\mn
As $p \le r_{\zeta(*)}$ this suffices. 

Toward contradiction assume that $r,\alpha$ are such that
\mn
\begin{enumerate}
\item[$\boxplus_3$]   $r_{\zeta(*)} \le r \in \Bbb P$ and 
$r \Vdash ``\alpha \in \bigcap\limits_{\zeta < \sigma^+} 
\name A_{\name \beta_\zeta}"$.
\end{enumerate}
\mn
Recall clause $(\delta)$ of $\boxplus_2$ (and $\zeta < \sigma^+ \Rightarrow
p_\zeta \le_{\text{\rm pr}} p_{\sigma^+} = q$), we know $\zeta <
\sigma^+ \Rightarrow \text{ Dom}(q) \cap (\beta^*_\zeta/E) \ne
\emptyset$ and, of course, $q \le r_{\zeta(*)} \le r$ so by
the definition of $\le$ in $\Bbb P$, for every $\xi < \sigma^+$ large enough
\mn
\begin{enumerate}
\item[$\boxplus_4$]   $(a) \quad (\beta^*_\xi/E) \cap \text{ Dom}(r_\xi) 
\backslash \text{Dom}(r_{\zeta(*)}) = \emptyset$ hence
\sn
\item[${{}}$]  $(b) \quad r,r_\xi$ are compatible functions (hence
conditions)
\sn
\item[${{}}$]  $(c) \quad \alpha \notin v_\zeta$.
\end{enumerate}
\mn
Let $r^+ = r \cup r_\zeta \cup \{\langle \kappa \beta^*_\zeta + \alpha,0
\rangle\}$.  

So easily
\mn
\begin{enumerate}
\item[$\boxplus_5$]  $(a) \quad r \le r^+ \in \Bbb P$
\sn
\item[${{}}$]  $(b) \quad r_\zeta \le r^+ \text{ hence } r^+ 
\Vdash "\name \beta_\zeta = \beta^*_\zeta"$
\sn
\item[${{}}$]  $(c) \quad r^+ \Vdash ``\alpha \notin 
\name A_{\name \beta^*_\zeta}"$.
\end{enumerate}
\mn
So we have gotten a contradiction thus proving $(*)_8$.
\mn
\begin{enumerate}
\item[$(*)_9$]   $\Vdash ``\neg\text{Pr}_2(\kappa^{++},\kappa,\theta,
\theta)$ if $\theta < \sigma"$.
\end{enumerate}
\mn
Why?  So toward contradiction suppose 
$p^* \Vdash ``\name {\cA} = \{\name B_\alpha:\alpha < \kappa^{++}\}
\subseteq [\kappa]^\kappa$ exemplifies 
Pr$_2(\kappa^{++},\kappa,\theta,\theta)"$.

For each $\alpha < \kappa^{++}$ we can find $p_\alpha$ such that $\bar
p,\bar r$ and then $S,q$:
\mn
\begin{enumerate}
\item[$\circledast_1$]  $(a) \quad \bar p  = \langle p_{\alpha,i}:i
\le \kappa\rangle$ is $\le_{\text{pr}}$-increasing continuous (in
$\bbP$)
\sn
\item[${{}}$]  $(b) \quad p_0 = p^*$
\sn
\item[${{}}$]  $(c) \quad p_{\alpha,i} \le r_{\alpha,i}$ and
$r_{\alpha,i} \Vdash_{\bbP} ``\gamma_{\alpha,i} \in \name B_\alpha
\backslash i"$
\sn
\item[${{}}$]  $(d) \quad p_{\alpha,i} \le_{\text{pr}} p_{\alpha,i+1}
\le_{\text{ap}} r_{\alpha,i}$
\sn
\item[${{}}$]  $(e) \quad S_\alpha \subseteq \kappa$ is stationary
\sn
\item[${{}}$]  $(f) \quad \langle v_{\alpha,i} := \text{
Dom}(r_{\alpha,i}) \backslash \text{ Dom}(p_{\alpha,i+1}):i \in
S\rangle$ is a $\Delta$-system with heart

\hskip25pt  $v_\alpha$, so $|v_\alpha| <
\sigma$
\sn
\item[${{}}$]  $(g) \quad \langle r_{\alpha,i} \rest v_\alpha:i \in
s\rangle$ is constantly $r_\alpha$
\sn
\item[${{}}$]  $(h) \quad p_\alpha = p_{\alpha,\kappa} \cup r_\alpha
\in \bbP$ so $p_\kappa \le q$.
\end{enumerate}
\mn
Let $u_\alpha = \cup\{\beta/E:\beta \in \text{ Dom}(p_\alpha)\}$ so
$u_\alpha \in [\chi]^{\le \kappa}$ and Dom$(r_{\alpha,i}) \subseteq
u_\alpha$ for $i < \kappa$.  For some 
$Y \in [\kappa^{++}]^{\kappa^{++}}$ and $\Upsilon^* < \kappa^+$ and
stationary $S \subseteq \kappa$ and $\bar\gamma = \langle \gamma_i:i
\in S\rangle$ we have
\mn
\begin{enumerate}
\item[$\circledast_2$]   if $\alpha \in Y$ then otp$(u_\alpha) =
\Upsilon^*$ and $S_\alpha = S$ and $\langle \gamma_{\alpha,i}:i \in
S\rangle = \bar\gamma$.
\end{enumerate}
\mn
Let $g_{\alpha,\beta}$ be the order preserving function from $u_\beta$ onto
$u_\alpha$.

Again as $2^\kappa = \kappa^+$ \wilog
\mn
\begin{enumerate}
\item[$\circledast_3$]  For $\alpha,\beta \in Y$
\begin{enumerate}
\item[$(a)$]  $p_\beta = p_\alpha \circ g_{\alpha,\beta}$ and
$p_{\beta,i} = p_{\alpha,i} \circ g_{\alpha,\beta}$ and $r_{\beta,i} =
r_{\alpha,i} \circ g_{\alpha,\beta}$ for $i < \kappa$
\sn
\item[$(b)$]  $u_\alpha \cap u_\beta = u_*$ for $\alpha
< \beta < \kappa^{++}$
\sn
\item[$(c)$]  $g_{\alpha,\beta}$ is the identity on $u_*$ 
\end{enumerate}
\end{enumerate}
and
\mn
\begin{enumerate}
\item[$\circledast_4$]   if $p_\alpha \le_{\text{apr}} r_\alpha,
p_\beta \le_{\text{apr}} r_\beta,r_\beta = r_\alpha \circ 
g_{\alpha,\beta}$ and $\gamma < \kappa$ then
\begin{enumerate}
\item[$(a)$]  $r_\alpha \Vdash 
``\gamma \in \name B_\alpha" \Leftrightarrow r_\beta \Vdash 
``\gamma \in \name B_\gamma"$
\sn
\item[$(b)$]  $r_\alpha \Vdash "\gamma \notin \name B_\alpha" \Leftrightarrow
r_\beta \Vdash ``\gamma \notin \name B_\gamma"$.
\end{enumerate}
\end{enumerate}
\mn
Choose $\langle \beta_\varepsilon:\varepsilon < \kappa^{++} \rangle$
an increasing sequence of ordinals from $Y$.

Let $p^* = \bigcup\limits_{\varepsilon < \theta} p_{\beta_\varepsilon}$ and
$\zeta(*) = \beta_\theta$.

So:
\mn
\begin{enumerate}
\item[$\circledast_5$]   $(a) \quad \Bbb P 
\models ``p^* \le p_{\beta_\varepsilon} < q^*"$ for $\varepsilon < \theta$
\sn
\item[${{}}$]  $(b) \quad$ if $p^* \le q$ then for each $\varepsilon <
\theta$ for every large enough $i \in S$, the 

\hskip25pt conditions $q,r_{\beta_\varepsilon,i}$ are compatible hence
\sn 
\item[${{}}$]  $(c) \quad$ if $q^* \le q$ then for every large enough
$i \in S$ the universal $q^+ =$

\hskip25pt $q \cup \{r_{\beta_\varepsilon,i}:
\varepsilon < \theta\}$ is a well
defined function, belongs to $\bbP$

\hskip25pt and is common upper bound of
$\{r_{\beta_\varepsilon,i}:\varepsilon  < \theta\} \cup q$ 

\hskip25pt hence force $\gamma_i \in B_{\beta_\varepsilon }$ 
for $\varepsilon  < \theta$.
\end{enumerate}
\mn
As $\gamma_i \ge i$ clearly
\mn
\begin{enumerate}
\item[$\circledast_6$]   $q^* \Vdash ``\bigcap\limits_{\varepsilon <
\theta} B_{\beta_\varepsilon}$ is unbounded in $\kappa$ hence has
cardinality $\kappa"$.
\end{enumerate}
\end{PROOF}

%\bibliographystyle{alphacolon}
%\bibliography{lista,listb,listx,listf,liste,listy}

\end{document}